\documentclass[a4paper,10pt,reqno]{amsart}

\usepackage[latin1]{inputenc}
\usepackage{amsmath, amssymb, amsthm, amsfonts, amsgen, bm, stmaryrd}
\numberwithin{equation}{section}
\usepackage{indentfirst}

\newcommand{\ds}{\displaystyle}

\newcommand{\Nb}{{\mathbb{N}}}
\newcommand{\Rb}{{\mathbb{R}}}
\newcommand{\Zb}{{\mathbb{Z}}}
\newcommand{\A}{{\mathcal{A}}}
\newcommand{\C}{{\mathcal{C}}}
\newcommand{\E}{{\mathcal{E}}}
\newcommand{\F}{{\mathcal{F}}}
\newcommand{\LL}{{\mathcal{L}}}
\newcommand{\M}{{\mathcal{M}}}
\newcommand{\PP}{{\mathcal{P}}}
\newcommand{\ie}{{\it i.e. }}
\newcommand{\med}{- \hskip -1em \int}

\makeatletter
\def\rightharpoonupfill@{\arrowfill@\relbar\relbar\rightharpoonup}
\newcommand{\xrightharpoonup}[2][]{\ext@arrow
0359\rightharpoonupfill@{#1}{#2}} \makeatother

\def\leq{\leqslant}
\def\geq{\geqslant}
\let\e=\varepsilon
\let\O=\Omega
\let\G=\Gamma
\let\a=\alpha
\let\b=\beta

\let\la=\langle
\let\ra=\rangle

\newtheorem{thm}{Theorem}[section]
\newtheorem{defi}[thm]{Definition}
\newtheorem{rmk}[thm]{Remark}
\newtheorem{lemma}[thm]{Lemma}
\newtheorem{proposition}[thm]{Proposition}
\newtheorem{coro}[thm]{Corollary}
\newtheorem{example}[thm]{Example}
\makeatletter

\begin{document}

\title[Characterization of two-scale gradient Young measures]
{Characterization of two-scale gradient Young measures and
application to homogenization}

\author[J.-F. Babadjian, M. Ba\'{\i}a \& P. M. Santos]
{Jean-Fran\c cois Babadjian, Margarida Ba\'{\i}a \& Pedro M. Santos}

\address[Jean-Fran\c{c}ois Babadjian]{SISSA, Via Beirut 2-4, 34014 Trieste, Italy}
\email{babadjia@sissa.it}

\address[Margarida Ba\'{\i}a]{Instituto Superior Tecnico, Av. Rovisco Pais, 1049-001 Lisbon, Portugal}
\email{mbaia@math.ist.utl.pt}

\address[Pedro M. Santos]{Instituto Superior Tecnico, Av. Rovisco Pais, 1049-001 Lisbon, Portugal}
\email{pmsantos@math.ist.utl.pt}

\maketitle

\begin{abstract}
This work is devoted to the study of two-scale gradient Young
measures naturally arising in nonlinear elasticity homogenization
problems. Precisely, a characterization of this class of measures
is derived and an integral representation formula for homogenized
energies, whose integrands satisfy very weak regularity
assumptions, is  obtained in terms of two-scale gradient Young
measures.
\end{abstract}

\begin{center}
\begin{minipage}{13cm}
\vspace{0.2cm}

\small{ 

\vspace{0.3cm}

 }
\end{minipage}
\end{center}

\vspace{0.5cm}

\section{Introduction}

\noindent Young (or Parametrized) measures have been introduced in
optimal control theory  by L. C. Young \cite{LCY3} to study non
convex variational problems for which there were no classical
solution, and to provide an effective notion of generalized solution
for problems in Calculus of Variations.

Starting with the works of Tartar \cite{T1} on hyperbolic
conservation laws, Young measures have been an important tool for
studying the asymptotic behavior of solutions of nonlinear partial
differential equations (see also DiPerna \cite{DP}). A key feature
of these measures is their capacity to capture the oscillations of
minimizing sequences of non convex variational problems, and many
applications arise {\it e.g.} in models of elastic crystals (see
Chipot \& Kinderlehrer \cite{CK} and Fonseca \cite{F}), phase
transition (see Ball \& James \cite{BJ}), optimal design (see
Bonnetier \& Conca \cite{BC},  Maestre \& Pedregal \cite{MP} and
Pedregal \cite{P3}). The special properties of Young measures
generated by sequences of gradients of Sobolev functions have been
studied by Kinderlehrer \& Pedregal \cite{KP2,KP}   and are relevant
in the applications to nonlinear elasticity.

The lack of information   on the spatial structure of oscillations
presents an obstacle for the application of Young measures to
homogenization problems. Two-scale Young measures, which have been
introduced by E in \cite{W} to study periodic homogenization of
nonlinear transport equations, contain some information on the
amount of oscillations and extend Nguetseng's notion of two-scale
convergence (see \cite{N1} and  Allaire \cite{A}). Other
(generalized) multiscale Young measures have been introduced in the
works of  Alberti \& M\"{u}ller \cite{AM} and Ambrosio $\&$ Frid
\cite{AF}.\\

From a variational point of view periodic homogenization of
integral functionals  rests on the study of the equilibrium
states, or minimizers, of a family of functionals of the type
\begin{equation}\label{main-funct}
\F_\e(u):=\int_\O f\left(x,\left\la \frac{x}{\e}\right\ra,\nabla
u(x)\right)\, dx,
\end{equation} as $\e\to 0$, under suitable boundary conditions. Here $\O$ (bounded
open subset of $\Rb^N$) is the reference configuration of a
nonlinear elastic body   with   periodic microstructure and whose
heterogeneities scale like a small parameter $\e>0$. The function
$u\in W^{1,p}(\O;\Rb^d)$ stands for a deformation and $f:\O \times Q
\times \Rb^{d \times N} \to [0,+\infty)$, with $Q:=(0,1)^N$, is the
stored energy density of this body that is assumed to satisfy
standard $p$-coercivity and $p$-growth conditions, with $p>1$. The
presence of the term $\la x/\e \ra$ (fractional part of the vector
$x/\e$ componentwise) takes into account the   periodic
microstructure of the body leading the integrand of
(\ref{main-funct}) to be periodic with respect to that variable. The
macroscopic (or averaged) description of this material may be
understood by the $\G$-limit of (\ref{main-funct}) with respect to
the weak $W^{1,p}(\O;\Rb^d)$-topology   (or, equivalently, with
respect to the  $L^{p}(\O;\Rb^d)$-topology if $\O$ is, for instance,
Lipschitz) and it has already been studied by many authors in the
Sobolev setting. Namely, under several regularity assumptions on $f$
it has been proved that for all $u\in W^{1,p}(\O;\Rb^{d})$,
\begin{equation}\label{fhom-2}
{\mathcal F}_{\rm hom}(u):=\G\text{-}\lim_{\e\to 0} \mathcal F_{\e}(u) =\int_\O f_{\rm
hom}(x,\nabla u(x))\, dx,\end{equation} \noindent  where for all
$(x,\xi) \in \O \times \Rb^{d \times N}$
\begin{eqnarray}\label{fhom}
&&\hspace{-1cm}f_{\rm hom}(x,\xi)\nonumber\\
&&\hspace{-0.5cm}= \lim_{T \to +\infty}\inf_{\phi} \left\{
\med_{(0,T)^N}f(x,y,\xi + \nabla \phi(y))\, dy : \phi \in
W^{1,p}_0((0,T)^N;\Rb^d)\right\}
\end{eqnarray}
(see Ba\'{\i}a \& Fonseca \cite{BF,BF1}, Braides \cite{B1}, Braides
\& Defranceschi \cite{BD}, Marcellini \cite{Ma} and M\"uller
\cite{Mu}).

We also refer   to Anza Hafsa, Mandallena \& Michaille \cite{HMM}
where a formula for the function $f_{\rm hom}$  has been given in
terms of gradient Young measures.  In the convex case,  Barchiesi
\cite{Bar} and Pedregal \cite{P} have derived the same $\G$-limit
result (\ref{fhom-2}) with Young measures techniques.   The main
contribution in \cite{Bar} is to weaken, as most as possible, the
regularity of  $f$ that is assumed to be an ``admissible integrand"
in the sense of Valadier \cite{V} (see Definition \ref{admint}
below).   Using the same kind of arguments, Pedregal has extended
this result to the nonconvex case in \cite{P2}.\\

We note that solutions of
$$\min_{u=u_0 \text{ on }\partial \O}  \int_\O f_{\rm
hom}(x,\nabla u(x))\, dx$$ \noindent  only give an average of the
oscillations  that minimizing sequences may develop.

From a mathematical
point of view, the main property of Young measures is their
capability of describing the asymptotic behavior of integrals of the
form $$\int_{\O} f(v_{\e}(x))\,dx,$$ where $f$ is some nonlinear
function and $\{v_{\e}\}$ is an oscillating sequence.    To address
the homogenization of  (\ref{main-funct}) we consider Young measures
generated by sequences of the type $\{(\la \cdot / \e\ra, \nabla
u_\e)\}$, which are, roughly speaking, what we will call {\it
two-scale gradient Young measures}. From a physical point of view,
we seek to capture microstructures -- due to finer and finer
oscillations of minimizing sequences that cannot reach an optimal
state -- at a given scale $\e$ (period of the material
heterogeneities). In this way, the minima of the limit problem
captures two kinds of oscillations of the minimizing sequences:
those due to the periodic heterogeneities of the material and those
due to a possible multi-well structure.

Our main result gives a complete algebraic
characterization of two-scale gradient Young
measures   (see Definition \ref{MGYM} below)  in the spirit of  Kinderlehrer \& Pedregal \cite{KP}.
We derive this characterization in terms of a Jensen's inequality with
test functions in the space  $\E_p$ of continuous functions
$f:{\overline Q} \times \Rb^{d \times N} \to \Rb$ such that the
limit $$\lim_{|\xi|\to +\infty}\frac{f(y,\xi)}{1+|\xi|^p}$$ exists
uniformly with respect to $y \in \overline Q$. Namely, we prove  the
following result.

\begin{thm}\label{bbs}Let $\O$ be a bounded open subset
of $\Rb^N$ with Lipschitz boundary and let $\nu \in L^\infty_w(\O
\times Q;\M(\Rb^{d \times N}))$ be such that $\nu_{(x,y)} \in
\PP(\Rb^{d \times N})$ for a.e.\ $(x,y) \in \O \times Q$. The family
$\{\nu_{(x,y)}\}_{(x,y) \in \O \times Q}$ is a {\rm two-scale
gradient Young measure} if and only if the three conditions below
hold:
\begin{itemize}
\item[i)] there exist  $u \in W^{1,p}(\Omega; \Rb^d)$ and $u_1 \in
L^{p}(\Omega; W_{per}^{1,p}(Q;\Rb^d))$ such that
\begin{equation}\label{i}\int_{\Rb^{d \times N}} \xi\, d\nu_{(x,y)}(\xi) =
\nabla u(x) + \nabla_y u_1(x,y) \quad \text{ for a.e. }(x,y) \in \O
\times Q;
\end{equation}
\item[ii)] for every $f\in \E_p$
\begin{equation}\label{ii}
\int_Q \int_{\Rb^{d \times N}} f(y,\xi) \, d\nu_{(x,y)}(\xi)\, dy
\geq f_{\rm hom}(\nabla u(x)) \quad \text{ for a.e. }x \in \O,
\end{equation}
where
\begin{eqnarray}\label{fhom1}
&&\hspace{-1cm}f_{\rm hom}(\xi)\nonumber\\
&&\hspace{-0.5cm}= \lim_{T \to +\infty}\inf_{\phi} \left\{
\med_{(0,T)^N}f(\la y \ra,\xi + \nabla \phi(y))\, dy : \phi \in
W^{1,p}_0((0,T)^N;\Rb^d)\right\};
\end{eqnarray}
\item[iii)] \begin{equation}\label{iii}\ds (x,y) \mapsto
\int_{\Rb^{d \times N}} {|\xi|}^p d\nu_{(x,y)}(\xi) \in L^1(\O
\times Q).\end{equation}
\end{itemize}
\end{thm}

We note that $\E_p$ is separable (see Section 3), and thus condition
(ii) needs only to be checked for countably many test functions $f$.
The proof of this theorem is similar to that of Kinderlehrer \&
Pedregal \cite{KP}. We first address the homogeneous case, that is,
we consider two-scale gradient Young measures that are independent
of the macroscopic variable $x \in \O$. This case rests on the
Hahn-Banach Separation Theorem. The general case will be obtained by
splitting $\O$ into suitable small subsets and by approximating
these measures by  two-scale Young measures that are piecewise
constant with respect to the variable $x \in \O$.

Theorem \ref{bbs} turns out to be useful to obtain a representation
of the $\G$-limit of (\ref{main-funct}) in terms of two-scale
gradient Young measures. This is the aim of our second result, where
following Barchiesi \cite{Bar}, we consider very weak regularity
hypothesis on the integrand $f$.

\begin{thm}\label{exprGlim}   Let $\O$ be a bounded open subset
of $\Rb^N$ with Lipschitz boundary and let  $f:\O\times Q\times
\Rb^{d\times N}\rightarrow [0,+\infty)$ be an admissible
integrand. Assume that there exist constants $\a$, $\b>0$ and $p
\in (1,+\infty)$ such that for all $(x,y,\xi) \in \O\times Q
\times \Rb^{d\times N}$
\begin{equation}\label{pgrowth}
\a|\xi|^p \leq f(x,y,\xi)\leq \b(1+|\xi|^{p}).
\end{equation}
Then  the functional  $\F_\e$ $\G$-converges with respect to the
weak $W^{1,p}(\O;\Rb^d)$-topology (or equivalently the strong
$L^p(\O;\Rb^d)$-topology) to $\F_{\rm hom}: W^{1,p}(\O;\Rb^{d})\to
[0,+\infty)$ given by
\begin{equation}\label{for-hom}\mathcal
F_{\rm hom}(u)= \min_{\nu \in \M_u} \int_\O \int_Q \int_{\Rb^{d
\times N}} f(x,y,\xi)\, d\nu_{(x,y)}(\xi)\,dy\, dx,\end{equation}
where
\begin{eqnarray}\label{mu}
{\mathcal M}_{u} & := & \Bigg\{ \nu \in L_w^\infty(\O \times Q ;
\M(\Rb^{d \times N})): ~ \{\nu_{(x,y)}\}_{(x,y)\in \O\times Q}
\text{ is a }\nonumber\\
&& \hspace{0.5cm}\text{two-scale gradient Young measure such that }\nonumber\\
&& \hspace{0.5cm}\nabla u(x) = \int_Q \int_{\Rb^{d \times N}}\xi
\, d\nu_{(x,y)}(\xi)\, dy ~\text{ for a.e. }x \in \O \Bigg\}
\end{eqnarray}
for all $u \in W^{1,p}(\O;\Rb^{d}).$
\end{thm}

The overall plan of this work in the ensuing sections will be as
follows: Section 2 collects the main notations and results used
throughout. Section 3 is devoted to the proof of Theorem \ref{bbs},
and in Section 4 we address the proof of the homogenization result
Theorem \ref{exprGlim}.

\section{Some preliminaries}

\noindent The purpose of this section is to give a brief overview of
the concepts and results that are used in the sequel. Almost all
these results are stated without proofs as they can be readily found
in the references given below.

\subsection{Notation}

\noindent Throughout this work $\O$ is an open bounded subset of
$\Rb^N$ with Lipschitz boundary, $\A(\O)$ denotes the family of all
open subsets of $\O$, $\LL^N$ is the Lebesgue measure in $\Rb^N$,
$\Rb^{d\times N}$ is identified with the set of real $d\times N$
matrices and $Q:=(0,1)^N$ is the unit cube in $\Rb^N$. The symbols
$\la \cdot \ra$ and $[\cdot]$ stand, respectively, for the
fractional and integer part of a number, or a vector, componentwise.
The Dirac mass at a point $a\in \Rb^m$ is denoted by $\delta_a$. The
symbol
$\med_A$ stands for the average $\text{\small$\LL^N(A)^{-1}$}\int_A$.\\

Let $U$ be an open subset of $\Rb^m$. Then:
\begin{itemize}
\item $\C_c(U)$ is the space of continuous functions $f:U\to \Rb$
with compact support.
\item $\C_0(U)$ is the closure of $\C_c(U)$ for the uniform
convergence; it coincides with the space of all continuous
functions $f:U \to \Rb$ such that, for every $\eta >0$, there
exists a compact set $K_\eta \subset U$ with $|f| < \eta$ on $U
\setminus K_\eta$.
\item $\M(U)$ is the space of real-valued Radon measures with
finite total variation. We recall that by the Riesz Representation
Theorem   $\M(U)$ can be identified with the dual space of $\C_0(U)$
through the duality $$\la \mu ,\phi \ra = \int_U \phi \, d\mu,
\qquad\mu \in \M(U), \quad \phi \in \C_0(U).$$
\item ${\mathcal P}(U)$ denotes  the space of probability
measures on $U$, \ie the space of all $\mu \in \M(U)$ such that $\mu
\geq 0$ and $\mu(U) =1$.
\item   $L^1(\O;\C_0(U))$ is the space of maps $\phi : \O \to
\C_0(U)$ such that
\begin{itemize}
\item[{\small i})] $\phi$ is strongly measurable, \ie there
exists a sequence of simple functions $s_n : \O \to \C_0(U)$ such
that $\|s_n(x) - \phi(x)\|_{\C_0(U)} \to 0$ for a.e. $x \in \O$;
\item[{\small ii})] $x \mapsto \|\phi(x)\|_{\C_0(U)} \in L^1(\O)$.
\end{itemize}
We recall that the linear space spanned by $\{\varphi \otimes \psi :
\; \varphi \in L^1(\O) \text{ and } \psi \in \C_{0}(U) \}$ is dense
in $L^{1}(\O;\C_{0}(U))$.
\item $L^{\infty}_{w}(\O;\M(U))$ is the space of maps $\nu:
\O \to \M(U)$ such that
\begin{itemize}
\item[{\small i)}] $\nu$ is weak* measurable, \ie $x
\mapsto \langle \nu_{x}, \phi \rangle$ is measurable for every $\phi
\in \C_0(U);$
\item[{\small ii)}] $x \mapsto \|\nu_{x}\|_{\M(U)}\in
L^{\infty}(\O).$
\end{itemize}
The space $L^{\infty}_{w}(\O;\M(U))$ can be identified with the dual
of $L^{1}(\O;\C_{0}(U))$ through the duality
$$\la \mu ,\phi \ra
= \int_\O \int_U \phi(x,\xi) \, d\mu_x(\xi)\, dx,  \qquad \mu \in
L^{\infty}_{w}(\O;\M(U)), \quad \phi \in L^1(\O;\C_0(U)),$$ where
$\phi(x,\xi):=\phi(x)(\xi)$ for all $(x,\xi)\in \O\times U$. Hence
it can endowed with the weak* topology (see {\it e.g.}\ Theorem 2.11
in M\'alek, Ne\v cas, Rokyta \& R\r u\v zi\v cka \cite{MNRR}).
\item The space $W_{\rm per}^{1,p}(Q;\Rb^d)$ stands for the
$W^{1,p}$-closure of all functions $f \in \C^1(\Rb^N,\Rb^d)$ which
are $Q$-periodic.
\end{itemize}

\subsection{Young measures}

\noindent We recall here the notion of Young  measure and some of
its basic properties. We refer the reader to {Braides}
\cite{BLNotes}, {M\"uller} \cite{M}, Pedregal \cite{Pbook},
Roub\'{\i}\v{c}ek \cite{TR}, Valadier \cite{V1} and references
therein for a detailed description on the subject.

\begin{defi}{\rm (Young measure)} Let  $\nu\in L^{\infty}_{w}(\O;\M(\Rb^m))$
and let $z_n:\Omega \rightarrow \Rb^m$  be a sequence of measurable
functions. The family of measures $\{\nu_x\}_{x \in \O}$ is said to
be the {\rm Young measure} generated by  $\{z_n\}$  provided
$\nu_{x}\in {\mathcal P}(\Rb^m)$ for a.e.\- $x\in \O$ and
$$\delta_{z_n}\xrightharpoonup {*}{} \nu~\text{in}~ L^{\infty}_{w}(\O;\M(\Rb^m)),$$
\ie for all $\psi\in L^{1}(\O;\C_{0}(\Rb^m))$
$$\lim_{n \to +\infty}\int_{\O}\psi(x,z_n(x))\,dx=\int_{\O}\int_{\Rb^m}\psi(x,\xi)\,d\nu_{x}(\xi)\,dx.$$
\end{defi}
The family $\{\nu_x\}_{x \in \O}$ is said to be a {\it homogeneous}
Young measure   if the map $x \mapsto \nu_x$ is independent of $x$.
In this case the family $\{\nu_x\}_{x \in \O}$ is identified with a
single element $\nu$ of $\M(\Rb^m)$.

The following result asserts the existence of Young measures (see
Ball \cite{ball}).

\begin{thm}\label{young}
Let $\{z_n\}$ be a sequence of measurable functions $z_n :\O \to
\Rb^m$. Then  there exist a subsequence $\{z_{n_k}\}$ and $\nu \in
L^{\infty}_{w}(\O; \M(\Rb^{m}))$ with $\nu_x \geq 0$ for a.e.\ $x
\in \O$, such that $\delta_{z_{n_k}}\xrightharpoonup {*}{} \nu$ in
$L^{\infty}_{w}(\O;\M(\Rb^m))$  and the following properties hold:
\begin{itemize}
\item[(i)] $\|\nu_x\|_{\M(\Rb^m)} = \nu_x(\Rb^m) \leq 1$ for
a.e.\ $x \in \O$;\\
\item[(ii)] if \, ${\rm dist}(z_{n_k},K)\to 0$ in measure
for some closed set $K\subset \Rb^{m}$, then
${\rm Supp}(\nu_x) \subset K$ for a.e.\ $x \in \O$;\\
\item[(iii)] $\|\nu_x\|_{\M(\Rb^m)} = 1$ if and only if there
exists a Borel function $g : \Rb^m \to [0,+\infty]$ such that
$$\lim_{|\xi|\to +\infty}g(\xi) = +\infty \quad \text{ and } \quad
\sup_{k \in \Nb} \int_\O g(z_{n_k}(x))\, dx <+\infty;$$
\item[(iv)] if $f:\O \times \Rb^m \to [0,+\infty]$ is a normal
integrand, then $$\liminf_{k \to +\infty} \int_\O f(x,z_{n_k}(x))\,
dx \geq \int_\O \int_{\Rb^m} f(x,\xi)\, d\nu_x(\xi)\, dx;$$
\item[(v)]if $(iii)$ holds and if $f:\O \times \Rb^m \to
[0,+\infty]$ is a Carath\'eodory integrand such that the sequence
$\{f(\cdot,z_{n_k})\}$ is equi-integrable then
$$\lim_{k \to +\infty} \int_\O f(x,z_{n_k}(x))\,dx=\int_\O \int_{\Rb^m} f(x,\xi)\, d\nu_x(\xi)\, dx.$$
\end{itemize}
\end{thm}

\subsection{Two-scale gradient Young measures}

\noindent As remarked by Pedregal \cite{P}, regular Young measures
do not always provide enough information on the oscillations of a
certain sequence $\{v_\e\}$. To better understand oscillations that
occur at a given length scale $\e$ we may study  the Young measure
generated by the pair $\{(\la \cdot/\e \ra,v_\e)\}$. In this paper
we are interested in the case where $v_\e=\nabla u_\e$, for some
sequence $\{u_\e\} \subset W^{1,p}(\O;\Rb^d)$, with $1<p<\infty$.

Let $\mu \in L_w^\infty(\O; \M(\Rb^N \times \Rb^{d \times N}))$ and
$\{u_\e\} \subset W^{1,p}(\O;\Rb^d)$ be such that the pair $\{(\la
\cdot/\e \ra,\nabla u_\e)\}$ generates the Young measure
$\{\mu_{x}\}_{x\in \O}$. By an application of the Generalized
Riemann-Lebesgue Lemma (see {\it e.g.}\ Lemma 5.2 in Allaire
\cite{A} or Theorem 3 in Lukkassen, Nguetseng \& Wall \cite{LNW})
the sequence $\{\la \cdot/\e\ra\}$ generates the homogeneous Young
measure $dy:=\LL^N \lfloor Q$ (restriction of the Lebesgue measure
to $Q$). Then by the Disintegration Theorem (see Valadier \cite{V2})
there exists a map $\nu \in L^\infty_w(\O \times Q;\M(\Rb^{d \times
N}))$ with $\nu_{(x,y)} \in \PP(\Rb^{d \times N})$ for a.e.\ $(x,y)
\in \O \times Q$ and such that $\mu_x = \nu_{(x,y)} \otimes dy$ for
a.e. $x\in \O$, \ie
$$\int_{\Rb^N \times \Rb^{d \times N}} \phi(y,\xi)\, d\mu_x(y,\xi)=\int_Q
\int_{\Rb^{d \times N}} \phi(y,\xi)\, d\nu_{(x,y)}(\xi)\, dy$$ for
every $\phi \in \C_0(\Rb^N \times \Rb^{d \times N})$.

The family $\{\nu_{(x,y)}\}_{(x,y) \in \O \times Q}$ is referred in
\cite{P} as the {\it two-scale (gradient) Young measure} associated
to the sequence $\{\nabla u_\e\}$ at scale $\e$. More precisely, we
give the following definition.

\begin{defi}\label{MGYM}
Let $\nu \in L^\infty_w(\O \times Q;\M(\Rb^{d \times N}))$. The
family $\{\nu_{(x,y)}\}_{(x,y) \in \O \times Q}$ is said to be a
{\rm two-scale gradient Young measure} if $\nu_{(x,y)} \in
\PP(\Rb^{d \times N})$ for a.e.\ $(x,y) \in \O \times Q$ and if for
every sequence $\{\e_n\} \to 0$  there exists a bounded
sequence $\{u_n\}$ in $W^{1,p}(\O;\Rb^d)$ such that
$\{(\la\cdot/\e_n \ra,\nabla u_n)\}$ generates the Young measure
$\{\nu_{(x,y)} \otimes dy\}_{x\in \O}$ \ie for every $z \in L^1(\O)$
and $\varphi \in \C_0(\Rb^N \times \Rb^{d \times N})$,
$$\lim_{n \to +\infty}\int_\O z(x)\,
\varphi\left(\left\la\frac{x}{\e_n}\right\ra,\nabla u_n(x)\right) dx
= \int_\O \int_Q \int_{\Rb^{d \times N}}z(x) \, \varphi(y,\xi)\,
d\nu_{(x,y)}(\xi)\, dy\, dx.$$ In this case
$\{\nu_{(x,y)}\}_{(x,y)\in \O\times Q}$ is also called the two-scale
Young measure associated to $\{\nabla u_n\}$.
\end{defi}

\begin{example}\label{ex}
{\rm Let $\{\e_n\} \to 0$, and let $u:\O\to\Rb^{d}$ and
$u_{1}:\O \times \Rb^N \to \Rb^d$ be smooth functions such that
$u_1(x,\cdot)$ is $Q$-periodic for all $x \in \O$. Define
$$u_n(x):=u(x)+\e_n u_{1}\left(x,\frac{x}{\e_n}\right).$$
The two-scale gradient Young measure $\{\nu_{(x,y)}\}_{(x,y)\in
\O\times Q}$ associated to $\{ \nabla u_n\}$ is given by
$$\nu_{(x,y)}:=\delta_{\nabla u(x) + \nabla_{y} u_{1}(x,y)} \text{ for
all } (x,y)\in \O\times Q.$$ Indeed, let us show that
$\{(\la\cdot/\e_n\ra,\nabla u_n)\}$ generates the Young measure
$\{\nu_{(x,y)}\otimes dy\}_{x \in \O}$. First we note that $\nabla
u_n(x) = \nabla u(x) + \e_n \nabla_x u_1(x,x/\e_n) + \nabla_y
u_1(x,x/\e_n)$ for every $x \in \O$. As $\e_n \nabla_x u_1(\cdot,
\cdot/\e_n) \to 0$ strongly in $L^p(\O;\Rb^{d \times N})$ because $x
\mapsto \nabla_x u_1(x,x/\e)$ is weakly convergent in $L^p(\O;\Rb^{d
\times N})$ (see {\it e.g.}\ Example 3 in Lukkassen, Nguetseng \&
Wall \cite{LNW}), then in particular
$$(\la \cdot / \e_n \ra ,\nabla u_n(\cdot)) - \left(\la \cdot / \e_n \ra , \nabla u(\cdot)
+ \nabla_y u_1(\cdot,\cdot/\e_n)\right) =(0,\e_n \nabla_x
u_1(\cdot,\cdot/\e_n)) \to 0$$ in measure. Thus from Lemma 6.3 in
Pedregal \cite{Pbook} the sequences $$\{(\la \cdot / \e_n \ra
,\nabla u_n)\} \text{ and } \{\left(\la \cdot / \e_n \ra , \nabla
u(\cdot) + \nabla_y u_1(\cdot,\cdot/\e_n)\right)\}$$ generate the
same Young measure. By Riemann-Lebesgue's Lemma we have for every
$\psi \in L^1(\O)$ and every $\varphi \in \C_0(\Rb^N \times \Rb^{d
\times N})$
\begin{eqnarray*}
\lim_{n \to +\infty}\int_\O \psi(x) \, \varphi\left(
\left\la\frac{x}{\e_n}\right\ra, \nabla u(x) + \nabla_y
u_1\left(x,\frac{x}{\e_n}\right) \right)\, dx\\
= \int_\O \int_Q \psi(x)\, \varphi(y, \nabla u(x) + \nabla_y
u_1(x,y) )\, dy \, dx \end{eqnarray*} which proves the claim. }
\end{example}

\begin{example}\label{ex1}
{\rm Let $\{\e_n\} \to 0$, and let $u:\O\to\Rb^{d}$ and
$u_{2}:\O \times \Rb^N \times \Rb^N \to \Rb^d$ be smooth functions
such that $u_2(x,\cdot,\cdot)$ is separately $Q$-periodic with
respect to its second and third variable, for all $x \in \O$. Define
$$v_n(x):=u(x)+\e_n^2 u_{2}\left(x,\frac{x}{\e_n},\frac{x}{\e_n^2}\right).$$
Arguing as previously, both sequences $$\{(\la \cdot / \e_n \ra
,\nabla v_n)\}\text{ and }\{\left(\la \cdot / \e_n \ra , \nabla
u(\cdot) + \nabla_z u_2(\cdot,\cdot/\e_n,\cdot/\e_n^2)\right)\}$$
generate the same Young measure. Using once more the
Riemann-Lebesgue Lemma we have that for every $\psi \in L^1(\O)$ and
every $\varphi \in \C_0(\Rb^N \times \Rb^{d \times N})$
\begin{eqnarray*} \lim_{n \to +\infty}\int_\O \psi(x) \,
\varphi\left( \left\la\frac{x}{\e_n}\right\ra, \nabla u(x) +
\nabla_z u_2\left(x,\frac{x}{\e_n},\frac{x}{\e_n^2}\right) \right)\, dx\\
= \int_\O \int_Q \int_Q \psi(x)\, \varphi(y, \nabla u(x) + \nabla_z
u_2(x,y,z) )\, dz\, dy \, dx. \end{eqnarray*} Hence, the two-scale
Young measure associated to $\{\nabla v_{n}\}$ is $$\nu_{(x,y)}:=
\int_Q \delta_{\nabla u(x) + \nabla_z u_{2} (x,y,z)} \,dz,$$ for
a.e.\ $(x,y) \in \O \times Q$, \ie
$$\la \nu_{(x,y)}, \phi \ra= \int_Q \phi(\nabla u(x) + \nabla_z u_{2}
(x,y,z))\,dz \quad \text{ for all }\phi \in \C_0(\Rb^{d \times
N}).$$ Note that in this example we do not get a Dirac mass because
there are oscillations occurring at different scales than $\e_n$,
namely at scale $\e_n^2$, that the two-scale Young measure misses
(see Valadier in \cite{V} for more details). }
\end{example}

\begin{rmk}\label{underlying}
{\rm Let $\{\e_n\}$, $\{u_n\}$ and $\nu$ be as in Definition
\ref{MGYM}. Since $\nabla u_n$ do not change if we add or remove a
constant, there is no loss of generality to assume that all the
functions $u_n$ have zero average. Moreover, let $\{u_{n_k}\}$ be a
subsequence of $\{u_n\}$.  Then there exists a subsequence
$\{u_{n_{k_{j}}}\}$ and $u \in W^{1,p}(\O;\Rb^d)$ (with zero
average) such that $u_{n_{k_{j}}} \rightharpoonup u$ in
$W^{1,p}(\O;\Rb^d)$ and
$$\nabla u(x) = \int_Q \int_{\Rb^{d \times N}} \xi\,
d\nu_{(x,y)}(\xi)\, dy \quad \text{ a.e.\ in }\O$$ (see {\it e.g.}\
the proof of Lemma \ref{nec} below).  It follows that $u$ is uniquely
defined because if $v$ is the weak $W^{1,p}(\O;\Rb^d)$-limit of
another subsequence of $u_{n_k}$ then
$$\nabla u(x) = \int_Q \int_{\Rb^{d \times N}} \xi\,
d\nu_{(x,y)}(\xi)\, dy=\nabla v(x) \quad \text{ a.e. in }\O,$$ which
implies that $u=v$ since they both have zero average. As a
consequence $u_n \rightharpoonup u$ in $W^{1,p}(\O;\Rb^d)$ and we
can   show in a similar way that $u$ is also independent of the
sequence $\{\e_n\}$. The function $u$ is called {\it the underlying
deformation} of $\{\nu_{(x,y)}\}_{(x,y) \in \O \times Q}$.}
\end{rmk}

In the following lemma, we show that there is no loss of generality
to assume that sequences of generators in Definition \ref{MGYM}
match the boundary condition of the underlying deformation.

\begin{lemma}\label{boundary}
Let $\{\e_n \} \to 0$ and $\{u_n\} \subset W^{1,p}(\Omega;\Rb^d)$ be
such that $u_n \rightharpoonup u$ in $W^{1,p}(\Omega;\Rb^d)$ for
some $u \in W^{1,p}(\Omega;\Rb^d)$. Suppose that $\{(\la
\cdot/\e_n\ra, \nabla u_n )\}$ generates the Young measure $\{
\nu_{(x,y)} \otimes dy \}_{x \in \O}$. Then there exists a sequence
$\{v_n\} \subset W^{1,p}(\Omega;\Rb^d)$ such that $v_n
\rightharpoonup u$ in $W^{1,p}(\Omega;\Rb^d)$, $v_n=u$ on a
neighborhood of $\partial \O$ and $\{(\la \cdot/\e_n\ra, \nabla v_n
)\}$ also generates $\{\nu_{(x,y)} \otimes dy \}_{x \in \O}$.
\end{lemma}

\begin{proof}
For any $k \in \Nb$ let $\O_k:=\left\{x\in \O:~ {\rm dist}(x,\Rb^N
\setminus \O)>1/k\right\}$ and let $\Phi_k \in \C^{\infty}_
{c}(\O;[0,1])$ be a cut-off function such that
$$\Phi_k:=\left\{\begin{array}{l}
1 \quad \text{if} ~ x \in \O_k,\\
0 \quad \text{if} ~ x \in   \O \setminus \overline
\O_{k+1}
\end{array}\right.$$
and $|\nabla \Phi_k| \leq C\,k$,  for some constant $C>0$. Let $w_{n,k} \in W^{1,p}(\O;\Rb^d)$ be
given by
$$w_{n,k}:=(1-\Phi_k)u+\Phi_k u_n,$$
\noindent from where
$$\nabla w_{n,k}=(1-\Phi_k)\nabla u+\Phi_k \nabla
u_n+(u_n-u)\otimes \nabla \Phi_k.$$ Since $u_n \to u$ strongly in
$L^p(\O;\Rb^d)$ then
\begin{equation}\label{bdy1}
\lim_{k \to +\infty} \lim_{n \to +\infty}
{\|w_{n,k}-u\|}_{L^p(\O;\Rb^d)}=0
\end{equation}
and, as a consequence of
$$\lim_{k \to +\infty} \lim_{n \to +\infty} {\|(u_n-u)\otimes \nabla
\Phi_k\|}_{L^p(\O;\Rb^{d \times N})}=0,$$ it follows that
\begin{equation}\label{bdy3}
\sup_{n,k \in \Nb}\|\nabla w_{n,k}\|_{L^p(\O;\Rb^{d \times N})}
<+\infty.
\end{equation}
Let $z$ and $\varphi$ be in a countable dense
subset of $L^1(\O)$ and $\C_{0}(\Rb^N \times \Rb^{d\times N}),$
respectively.  Then
\begin{eqnarray}\label{bdy4}
&&\lim_{k \to +\infty} \lim_{n \to +\infty}\int_{\O}z(x)\,
\varphi\left(\left\la\frac{x}{\e_n}\right\ra, \nabla
w_{n,k}(x)\right)\, dx\nonumber\\
&&\hspace{2cm} = \lim_{k \to +\infty} \lim_{n \to +\infty}\int_{\O_k}z(x)\,
\varphi\left(\left\la\frac{x}{\e_n}\right\ra, \nabla
u_n(x)\right)\, dx\nonumber\\
&&\hspace{3cm} + \lim_{k \to +\infty} \lim_{n \to +\infty}\int_{\O
\setminus \O_k}z(x)\, \varphi\left(\left\la\frac{x}{\e_n}\right\ra,
\nabla w_{n,k}(x)\right)\, dx\nonumber\\
&&\hspace{2cm} = \int_{\O} z(x) \int_Q \int_{\Rb^{d \times N}}
\varphi (y,\xi) \, d\nu_{(x,y)} \, dy \, dx.
\end{eqnarray}
By a diagonalization argument (see {\it e.g.} Lemma 7.2 in Braides,
Fonseca \& Francfort \cite{BFF}) and taking into account
(\ref{bdy1})-(\ref{bdy4}), we can find a sequence $\{k(n)\} \nearrow
+\infty$ such that, upon setting $v_n:=w_{n,k(n)}$, we have $v_n
\rightharpoonup u$ in $W^{1,p}(\O;\Rb^d)$, $v_n=u$ on a neighborhood
of $\partial \O$ and for every $z$ and $\varphi$ in a countable
dense subset of $L^1(\O)$ and $\C_0(\Rb^N \times \Rb^{d\times N})$,
respectively,
\begin{equation*}
\lim_{n \to +\infty}\int_{\O}z(x)\,
\varphi\left(\left\la\frac{x}{\e_n}\right\ra, \nabla v_n(x)\right)\,
dx = \int_\Omega z(x) \int_Q \int_{\Rb^{d \times N}} \varphi(y,\xi)
\,d\nu_{(x,y)}(\xi) \,dy\,dx.
\end{equation*}
\end{proof}

A two-scale gradient Young measure $\{\nu_{(x,y)}\}_{(x,y) \in \O
\times Q}$ is said to be {\it homogeneous} if the map $(x,y) \mapsto
\nu_{(x,y)}$ is independent of $x$. In this case, $\nu$ can be
identified with an element of $L^\infty_w(Q;\M(\Rb^{d \times N}))$
and we write $\{\nu_y\}_{y \in Q} \equiv\{\nu_{(x,y)}\}_{(x,y) \in
\O \times Q}.$\\

We next define  the average  of a map $\nu \in L^\infty_w(\O \times
Q;\M(\Rb^{d \times N}))$ for which $\{\nu_{(x,y)}\}_{(x,y) \in \O
\times Q}$ is a two-scale gradient Young measure. This notion,
useful for the analysis developed on Section \ref{homogeneous}, will
provide an important example of homogeneous two-scale gradient Young
measure.

\begin{defi}\label{ave} Let $\nu \in L^\infty_w(\O \times Q;\M(\Rb^{d \times
N}))$ be such that $\{\nu_{(x,y)}\}_{(x,y) \in \O \times Q}$ is a
two-scale gradient Young measure. The average of
$\{\nu_{(x,y)}\}_{(x,y) \in \O \times Q}$ (with respect to the
variable $x$) is the   family  $\{\overline \nu_{y}\}_{y \in Q}$
defined by
$$\la {\overline \nu}_y, \varphi \ra:=
\med_{\O}\int_{\Rb^{d \times N}}\varphi(\xi)\,d\nu_{(x,y)} \,dx$$ for every
$\varphi \in \C_0(\Rb^{d \times N})$.
\end{defi}

If  $\{\nu_{(x,y)}\}_{(x,y) \in \O \times Q}$  is a two-scale
gradient Young measure, then it can be seen that   $\overline
\mu:= \overline \nu_y \otimes dy$ is the average of $\{\mu_x\}_{x
\in \O}$ with $\mu_x:=\nu_{(x,y)} \otimes dy$ and $\mu\in
L^{\infty}_{w}(\O;\M(\Rb^N \times \Rb^{d \times N}))$. Thus,
$\overline \mu$ is a homogeneous Young measure by Definition
\ref{MGYM} and Theorem 7.1 in Pedregal \cite{Pbook}.

In the following Lemma we prove that $\{\overline \nu_y\}_{y \in
Q}$ is actually a homogeneous two-scale gradient Young measure. We
will use the same kind of blow up argument as in the proof of
Theorem 7.1 in Pedregal \cite{Pbook}, splitting $Q$ into suitable
subsets. However, contrary to \cite{Pbook} we will not use
Vitali's Covering Theorem because the  radii  of this sets
(which may vary from   one to another) may interact with the
length scale of our problem, $\e$,  in a way we are unable to
control. We will construct a covering   consisting  of
subsets of fixed radius. It is enough for our purposes to consider
the case where the underlying deformation is affine and $\O=Q$.

\begin{lemma}\label{average}
Let $\nu \in L^\infty_w(Q \times Q;\M(\Rb^{d \times N}))$ be such
that $\{\nu_{(x,y)}\}_{(x,y) \in Q \times Q}$ is a two-scale
gradient Young measure with underlying deformation $F\cdot$, for
$F\in \Rb^{d\times N}$. Then $\{\overline \nu_y\}_{y\in Q}$ is a
homogeneous two-scale gradient Young measure with the same underlying
deformation.
\end{lemma}

\begin{proof}
Note that by Definition \ref{ave} and Fubini's Theorem, it follows
that $y \mapsto \overline \nu_y$ is weakly*-measurable and thus
$\overline \nu \in L^\infty_w(Q;\M(\Rb^{d \times N}))$. We have to
show that for every sequence $\{\e_n \} \to 0$, there exists
$\{v_n\} \subset W^{1,p}(Q;\Rb^{d})$ such that $\{(\la \cdot/\e_n
\ra,\nabla v_n)\}$ generates $\overline
\mu:=\overline{\nu}_{y}\otimes dy$ and $v_n \rightharpoonup F\cdot$
in $W^{1,p}(Q;\Rb^{d}).$

Let $\{u_n\} \subset W^{1,p}(Q;\Rb^d)$ be such that $\{(\la n \,
\cdot \ra, \nabla u_n)\}$ generates $\{\nu_{(x,y)} \otimes dy\}_{x
\in \O}$ and $u_n \rightharpoonup F\cdot$ in $W^{1,p}(Q;\Rb^d)$ (see
Remark \ref{underlying}). Without loss of generality we may assume
that   $u_n(x)= Fx$ on a neighborhood of $\partial Q$ (see Lemma
\ref{boundary}).

Let $\{\e_n\} \to 0$ and for each $n$ define
$\rho_n:=\e_n [1/\sqrt\e_n]$. Then there exist $m_n \in \Nb$, $a_i^n
\in \rho_n \Zb^N \cap Q$ and a measurable set $E_n \subset Q$ with
$\LL^N(E_n) \to 0$ such that
$$Q= \bigcup_{i=1}^{m_n}\big(a^n_{i}+\rho_{n}Q\big) \cup E_n.$$
Define
$$v_n(x):=\left\{ \begin{array}{ll}\rho_n
u_{\rho_n/\e_n}\left(\frac{x-a^n_{i}}{\rho_n}\right)+Fa^n_{i}&
\text{ if } x\in a^n_{i}+\rho_n Q \text{ and } i \in
\{1,\ldots,m_n\},\\\\
Fx  & \text{ otherwise}. \end{array}\right.$$ Note that the previous
definition makes sense since $\rho_n/\e_n \in \Nb$. We remark that
$\{v_n\} \subset W^{1,p}(Q;\Rb^d)$ and  $v_n \rightharpoonup F\cdot$
in $W^{1,p}(Q;\Rb^d)$ since $u_n  \rightharpoonup F\cdot$ in this
space. Let $z \in \C_c(Q)$ and $\varphi\in \C_{0}(\Rb^N \times
\Rb^{d\times N})$. Then we have
\begin{eqnarray*}
&&\int_{Q}z(x)\, \varphi\left(\left\la\frac{x}{\e_n}\right\ra,
\nabla v_n(x)\right)\, dx\\
&&\hspace{1cm}=  \sum_{i=1}^{m_n} \int_{a^n_{i}+\rho_n Q}z(x)\,
\varphi\left(\left\la\frac{x}{\e_n}\right\ra,
\nabla u_{\rho_n/\e_n}\left(\frac{x-a^n_{i}}{\rho_n}\right)\right)\, dx\\
&&\hspace{1.5cm} + \int_{E_n}z(x)\,
\varphi\left(\left\la\frac{x}{\e_n}\right\ra,
F\right)\, dx \\
&&\hspace{1cm}=\sum_{i=1}^{m_n}
z(a^n_{i})\int_{a^n_{i}+\rho_nQ}\varphi\left(\left\la\frac{x}{\e_n}\right\ra,
\nabla u_{\rho_n/\e_n}\left(\frac{x-a^n_{i}}{\rho_n}\right)\right)\,
dx\\
&&\hspace{1.5cm} +\sum_{i=1}^{m_n} \int_{a^n_{i}+\rho_nQ}
\left(z(x)-z(a^n_{i})\right)
\varphi\left(\left\la\frac{x}{\e_n}\right\ra, \nabla
u_{\rho_n/\e_n}\left(\frac{x-a^n_{i}}{\rho_n}\right)\right)\,
dx\\
&&\hspace{1.5cm} + \int_{E_n}z(x)\,
\varphi\left(\left\la\frac{x}{\e_n}\right\ra, F\right)\, dx,
\end{eqnarray*}
and, consequently,
\begin{eqnarray}\label{1333}
&&\int_{Q}z(x)\, \varphi\left(\left\la\frac{x}{\e_n}\right\ra,
\nabla v_n(x)\right)\, dx\nonumber\\
&&\hspace{2cm}= \sum_{i=1}^{m_n}\rho_n^N
z(a^n_{i})\int_{Q}\varphi\left(\left\la\frac{a^n_{i}+\rho_nx}{\e_n}\right\ra,
\nabla u_{\rho_n/\e_n}(x)\right)\, dx\nonumber\\
&&\hspace{3cm} + o(1), \text{ as }n \to +\infty
\end{eqnarray}
by changing variables, using the uniform continuity of $z$ and the
fact that $\LL^N(E_n) \to 0$. Hence, as $a_i^n/\e_n \in \Zb^N$, it
follows that
\begin{eqnarray}\label{1333}
&&\int_{Q}z(x)\, \varphi\left(\left\la\frac{x}{\e_n}\right\ra,
\nabla v_n(x)\right)\, dx\nonumber\\
&&\hspace{2cm}= \sum_{i=1}^{m_n}\rho_n^N
z(a^n_{i})\int_{Q}\varphi\left(\left\la\frac{x}{\e_n/\rho_n}\right\ra,
\nabla u_{\rho_n/\e_n}(x)\right)\, dx\nonumber\\
&&\hspace{3cm} + o(1), \text{ as }n \to +\infty,
\end{eqnarray}
and passing to the limit in (\ref{1333}) and using Definition
\ref{ave}, we conclude that
\begin{eqnarray*}
&&\lim_{n \to +\infty}\int_{Q}z(x)\,
\varphi\left(\left\la\frac{x}{\e_n}\right\ra, \nabla v_n(x)\right)\,
dx\\
&&\hspace{2cm}= \int_{Q} z(x)\, dx \int_{Q} \int_Q \int_{\Rb^{d
\times N}}
\varphi(y,\xi)\, d\nu_{(x,y)}(\xi)\,dy \, dx\\
&&\hspace{2cm}= \la \overline \nu, \varphi \ra \int_{Q} z(x)\, dx.
\end{eqnarray*}
Since by density the previous equality holds for every $z \in
L^1(\O)$, then  $\{(\la\cdot/\e_n\ra,\nabla v_n)\}$ generates the
homogeneous Young measure $\overline \nu_y \otimes dy$ and,
consequently, $\{\overline\nu_y\}_{y \in Q}$ is a homogeneous
two-scale gradient Young measure.
\end{proof}

\section{Proof of Theorem \ref{bbs}}\label{Proof1}

\noindent The aim of this section is to prove Theorem \ref{bbs}. We
start by  introducing the space $\E_p$ of all continuous functions
$f : \overline Q \times \Rb^{d \times N} \to \Rb$ such that the
limit $$\lim_{|\xi| \to +\infty} \frac {f(y,\xi)}{1+{|\xi|}^p}$$
exists uniformly with respect to $y\in \overline Q$. As an example,
the function $(y,\xi) \mapsto a(y)|\xi|^p$, where $a \in
\C(\overline Q)$, is in $\E_p$.

It can be checked that $\E_p$ is a Banach space under the norm
\begin{equation*}
{\|f\|}_{\E_p}:= \sup_{y \in \overline Q, \, \xi \in \Rb^{d \times
N}} \frac {|f(y,\xi)|}{1+{|\xi|}^p}.
\end{equation*}
In addition, $\E_p$ is isomorphic to the space $\C\big({\overline Q}
\times (\Rb^{d \times N} \cup \{\infty\})\big)$ under the map
$$\begin{array}{llll}
 & \E_p & \longrightarrow &  \C\big({\overline Q} \times
(\Rb^{d \times N} \cup \{\infty\}) \big)\\
\ds & f & \longmapsto & \ds (y,\xi) \mapsto
\left\{\begin{array}{lll}
\hspace{0.3cm} \frac{f(y,\xi)}{1+|\xi|^p}& \text{ if }&(y,\xi)\in \overline{Q}\times \Rb^{d\times N}\\
\lim\limits_{|\xi| \to +\infty}  \frac {f(y,\xi)}{1+{|\xi|}^p}&
\text{ if } & |\xi|=+\infty,
\end{array}\right.
\end{array}$$
where $\Rb^{d \times N} \cup \{\infty\}$ denotes the one-point
compactification of $\Rb^{d \times N}$, and, consequently, it is
separable. Furthermore, for all $f\in \E_p$ there exists a constant
$c>0$ such that
\begin{equation}\label{Ep}
|f(y,\xi)| \leq c(1+|\xi|^p), \quad \text{ for all }(y,\xi) \in
\overline Q \times \Rb^{d \times N}.
\end{equation}

We denote by $(\E_p)'$ the dual space of $\E_p$ and the brackets
$\la \cdot,\cdot\ra_{(\E_p)',\E_p}$ stand for the duality product
between $(\E_p)'$ and $\E_p$.

\subsection{Necessity}

\noindent We start by showing that conditions i)-iii) in
(\ref{i})-(\ref{iii}) are necessary. Precisely we prove the
following result.

\begin{lemma}\label{nec} Let $\nu \in L^\infty_w(\O \times Q;\M(\Rb^{d \times
N}))$ be such that  $\{\nu_{(x,y)}\}_{(x,y) \in \O \times Q}$ is a
{\rm two-scale gradient Young measure}. Then
\begin{itemize}
\item[i)] there exist  $u \in W^{1,p}(\Omega; \Rb^d)$ and $u_1 \in
L^{p}(\Omega; W_{per}^{1,p}(Q;\Rb^d))$ such that
$$\int_{\Rb^{d \times N}} \xi \, d\nu_{(x,y)}(\xi) = \nabla u(x) + \nabla_y u_1(x,y)
 \quad \text{ for a.e.\ }(x,y) \in \O \times Q;$$
\item[ii)] for every  $f\in \E_p$  we have that
\begin{equation*}
\int_Q \int_{\Rb^{d \times N}} f(y,\xi) \, d\nu_{(x,y)}(\xi)\, dy
\geq f_{\rm hom}(\nabla u(x)) \quad \text{ for a.e.\ }x \in \O,
\end{equation*}
where $f_{\rm hom}$ is given by (\ref{fhom1}); \item[iii)] $\ds
(x,y) \mapsto \int_{\Rb^{d \times N}} {|\xi|}^p d\nu_{(x,y)}(\xi)
\in L^1(\O \times Q).$
\end{itemize}
\end{lemma}

\begin{proof}
Let $\{\nu_{(x,y)}\}_{(x,y) \in\O \times Q}$ be a two-scale gradient
Young measure.

We start by proving  that i) holds. By Definition \ref{MGYM} and
Remark \ref{underlying} there exists $u \in W^{1,p}(\O;\Rb^d)$ such
that for every sequence $\{\e_n\} \to 0$ one can find $\{u_n\}
\subset W^{1,p}(\O;\Rb^d)$ such that $\{(\la \cdot/\e_n\ra ,\nabla
u_n) \}$ generates the Young measure $\{\nu_{(x,y)} \otimes dy\}_{x
\in \O}$ and $u_n \rightharpoonup u$ in $W^{1,p}(\O;\Rb^d)$. Up to a
subsequence (still denoted by $u_n$), we can also assume that
$\{|\nabla u_n|^p\}$ is equi-integrable (see the Decomposition Lemma
in Fonseca, M\"uller \& Pedregal \cite{FMP}) and that there exists a
function $u_1 \in L^{p}(\O;W^{1,p}_{\rm per} (Q;\Rb^d))$ such that
the sequence $\{\nabla u_n\}$ two-scale converges to $\nabla u +
\nabla_y u_1$ (see {\it e.g.}\ Theorem 13 in Lukkassen, Nguetseng \&
Wall \cite{LNW}; see also Allaire \cite{A} or Nguetseng \cite{N1}).
Consequently, for all $\phi \in \C^\infty_c(\O \times Q;\Rb^{d
\times N})$ we have that
\begin{eqnarray}\label{2scaleconv}&&\lim_{n \to +\infty} \int_\O \nabla u_n(x) \cdot
\phi\left(x,\left\la\frac{x}{\e_n}\right\ra\right)dx\nonumber\\
&&\hspace{2cm} = \int_\O \int_Q (\nabla u(x) + \nabla_y u_1(x,y))
\cdot \phi(x,y)\, dy\, dx.\end{eqnarray} Set $f(x,y,\xi)=\xi \cdot
\phi(x,y)\,\,$ for $(x,y,\xi) \in \O\times Q\times \Rb^{d \times
N}$. As $f$ is a Carath\'eodory integrand (measurable in $x$ and
continuous in $(y,\xi)$) and the sequence $\{f(\cdot,\la
\cdot/\e_n\ra,\nabla u_n(\cdot))\}$ is equi-integrable, by Theorem
\ref{young} v) we get that
\begin{equation}\label{YM} \lim_{n \to +\infty} \int_\O \nabla u_n(x) \cdot
\phi\left(x,\left\la\frac{x}{\e_n}\right\ra\right)dx= \int_\O \int_Q
\int_{\Rb^{d \times N}} \xi \cdot \phi(x,y)\, d\nu_{(x,y)}(\xi)\,
dy\, dx.
\end{equation}
Consequently, from (\ref{2scaleconv}) and (\ref{YM}) we get for
a.e.\ $(x,y) \in \O \times Q$
$$\int_{\Rb^{d \times N}} \xi \, d\nu_{(x,y)}(\xi) = \nabla u(x) +
\nabla_y u_1(x,y),$$ which proves i).

Let us see now that iii) is satisfied. As  $\{\nabla u_n\}$ is
$p$-equi-integrable then by Theorem \ref{young} v) we get that
$$\int_\O \int_Q \int_{\Rb^{d \times N}} |\xi|^p \,
d\nu_{(x,y)}(\xi)\, dy\, dx = \lim_{n \to +\infty} \int_\O |\nabla
u_n|^p\, dx < +\infty,$$ which completes the proof of  iii).

Finally, let us see that condition ii) holds by application of the
classical $\G$-convergence result for the homogenization of integral
functionals (see Braides \cite{B1} or M\"{u}ller \cite{Mu}). Let $f
\in \E_p$. In particular $f$ satisfies the $p$-growth condition
(\ref{Ep}) but it is not necessarily $p$-coercive. For every $\a>0$
and $M>0$, define $f_{M,\a}(y,\xi):=f_M(y,\xi) + \a |\xi|^p$ where
$f_M(y,\xi)=\max\{-M,f(y,\xi)\}$. Then
$$\a|\xi|^p - M \leq
f_{M,\a}(y,\xi) \leq (c+\a) (1+|\xi|^p), \quad \text{ for all
}(y,\xi) \in \overline Q \times \Rb^{d \times N}.$$ Hence, by e.g.\
Theorem 14.5 in Braides \cite{B1}   ($\G$-$\liminf$ inequality) and
since $f_{M,\a} \geq f$, we get that for every $A \in \A(\O)$
\begin{eqnarray}\label{2056}
\liminf_{n \to +\infty} \int_A f_{M,\a} \left(\left\la
\frac{x}{\e_n}\right\ra,\nabla u_n(x) \right)dx & \geq & \int_A
(f_{M,\a})_{\rm hom}(\nabla u(x))\, dx\nonumber\\
& \geq & \int_A f_{\rm hom}(\nabla u(x))\, dx
\end{eqnarray}
where $f_{\rm hom}$ is defined in (\ref{fhom1}). On the other hand,
\begin{eqnarray}\label{2057}
\liminf_{n \to +\infty} \int_A f_{M,\a}\left(\left\la
\frac{x}{\e_n}\right\ra,\nabla u_n(x) \right)dx  & \leq & \liminf_{n
\to +\infty} \int_A f_M \left(\left\la
\frac{x}{\e_n}\right\ra,\nabla u_n(x) \right)dx\nonumber\\
&& + \a \sup_{n \in \Nb} \int_A |\nabla u_n|^p\, dx.
\end{eqnarray}
Gathering (\ref{2056}) and (\ref{2057}), and passing to the limit
as $\a \to 0$, we obtain that
\begin{equation}\label{gliminfca}
\liminf_{n \to +\infty} \int_A f_M\left(\left\la
\frac{x}{\e_n}\right\ra,\nabla u_n (x)\right)dx \geq \int_A f_{\rm
hom}(\nabla u(x))\, dx.
\end{equation}
Define the set
$$A_n^M:= \left\{x \in A : \; f\left(\left\la
\frac{x}{\e_n}\right\ra, \nabla u_n(x) \right) \leq -M \right\}$$
and notice that by Chebyshev's Inequality  $$\LL^N (A_n^M) \leq
c/M,$$ for some constant $c >0$ independent of $n$ and $M$. Then
\begin{eqnarray}\label{cali}
\int_A f_M \left(\left\la \frac{x}{\e_n}\right\ra, \nabla u_n(x)
\right) dx & = & -M \LL^N(A_n^M)\nonumber\\
&& + \int_{A\setminus A_n^M} f\left(\left\la
\frac{x}{\e_n}\right\ra, \nabla u_n(x) \right)dx\nonumber\\
& \leq &  \int_{A\setminus A_n^M} f\left(\left\la
\frac{x}{\e_n}\right\ra, \nabla u_n(x) \right)dx.
\end{eqnarray}
As $\{ |\nabla u_n|^p \}$ is equi-integrable, by the $p$-growth
condition (\ref{Ep}), it follows that $\{f(\la \cdot/\e_n\ra,\nabla
u_n)\}$ is also equi-integrable. Thus
\begin{equation}\label{matin}
\int_{A_n^M}f\left(\left\la \frac{x}{\e_n}\right\ra, \nabla u_n(x)
\right)dx \xrightarrow[M \to +\infty]{} 0 \end{equation} uniformly
with respect to $n \in \Nb$. By (\ref{gliminfca}), (\ref{cali})
and (\ref{matin}) we get that
\begin{equation}\label{gliminfbis}
\liminf_{n \to +\infty} \int_A f\left(\left\la
\frac{x}{\e_n}\right\ra,\nabla u_n (x)\right)dx \geq \int_A f_{\rm
hom}(\nabla u(x))\, dx.
\end{equation}
Finally, since $\{f(\la \cdot/\e_n\ra,\nabla u_n)\}$ is
equi-integrable,  by Theorem \ref{young} v) we have that
\begin{equation}\label{2scale}
\lim_{n \to +\infty} \int_A
f\left(\left\la\frac{x}{\e_n}\right\ra,\nabla u_n (x)\right)dx =
\int_A \int_Q \int_{\Rb^{d \times N}} f(y,\xi)\, d\nu_{(x,y)}(\xi)\,
dy\, dx
\end{equation}
and we conclude the proof of ii) thanks to (\ref{gliminfbis}) and
(\ref{2scale}) together with a localization argument.
\end{proof}

\subsection{Sufficiency}

\noindent We show here that these conditions are also sufficient to
characterize two-scale gradient Young measures. Following the lines
of Kinderlehrer \& Pedregal \cite{KP}, we first study the
homogeneous case. The non-homogeneous one will be obtained through a
suitable approximation of two-scale gradient Young measures by
piecewise constant ones.

\subsubsection{Homogeneous case}\label{homogeneous}

\noindent Our aim here is to prove the following result.
\begin{lemma}\label{homo}
Let $F \in \Rb^{d \times N}$ and $\nu \in L_w^\infty(Q;\M(\Rb^{d
\times N}))$ be such that $\nu_y \in \PP(\Rb^{d \times N})$ for
a.e.\ $y\in Q$. Assume that
\begin{equation}\label{n1}\displaystyle F =
\int_Q \int_{\Rb^{d \times N}}\xi \, d\nu_y(\xi)\,
dy,\end{equation}
 \begin{equation}\label{n2}f_{\rm hom}(F) \leq \int_Q
\int_{\Rb^{d \times N}} f(y,\xi)\, d\nu_y(\xi)\; dy \end{equation}
for every $f \in \E_p$, and that
\begin{equation}\label{powerp}\int_Q \int_{\Rb^{d \times N}}
|\xi|^p \, d\nu_y(\xi)\, dy <+\infty.\end{equation} Then
$\{\nu_y\}_{y \in Q}$ is a homogeneous two-scale  gradient Young
measure.
\end{lemma}

As Kinderlehrer \&  Pedregal \cite{KP}, the argument in this case
will rest on the Hahn-Banach Separation Theorem that implies any
element $\nu \in L^\infty_w(Q;\M(\Rb^{d \times N}))$, for which the
hypothesis of Theorem \ref{bbs} are satisfied, to be in a suitable
convex and weak* closed subset of homogeneous two-scale gradient
Young measures. To prove Lemma \ref{homo} we start by giving some
notations and auxiliary lemmas.

For every $F \in \Rb^{d \times N}$ let
\begin{eqnarray}\label{MF}
M_F & \hspace{-0.8cm} := \hspace{-0.6cm}& \bigg\{\nu \in
L^\infty_w(Q;\M(\Rb^{d \times
N})): ~~~~\{\nu_y\}_{y \in Q} \text{ is a homogeneous two-scale }\nonumber\\
&&\hspace{0.5cm}\text{gradient Young measure and } \int_Q
\int_{\Rb^{d \times N}}\xi\, d\nu_{y}(\xi)\, dy = F \bigg\}.
\end{eqnarray}

\begin{rmk}\label{indep}{\rm
The set $M_F$ is independent of $\O$, \ie if $\nu \in M_F$ and
$\O'\subset \Rb^{N}$ is another domain, then for all $\{\e_n\} \to
0$ there exist a sequence $\{v_n\}\in W^{1,p}(\O';\Rb^d)$ such that
$\{( \la\cdot / \e_n \ra,\nabla v_n)\}$ generates $\nu_y \otimes
dy$. Indeed, let   $r>0$ such that $\O' \subset r \O$. Fix an
arbitrary sequence $\{\e_n\} \to 0$ and define $\delta_n=\e_n/r$.
Then there exists a sequence $\{u_n\}\subset W^{1,p}(\O;\Rb^d)$ such
that $\{( \la \cdot / \delta_n \ra,\nabla u_n)\}$ generates the
homogeneous Young measure $\nu_y \otimes dy$. Define now $v_n(x)=r
\, u_n(x/r)$ so that $v_n$ belongs to $W^{1,p}(r\O ;\Rb^d)$ and thus
{\it a fortiori} to $W^{1,p}(\O';\Rb^d)$. A simple change of
variable shows that the sequence $\{(\la \cdot / \e_n \ra, \nabla
v_n)\}$ generates the homogeneous Young measure $\nu_y \otimes dy$
as well. }\end{rmk}

The next technical result allows us to construct two-scale
gradient Young measures from measures of this class that are
defined on disjoint subsets of $\O$. It will be of use in Lemma
\ref{closeconvex} below to prove the convexity of the set $M_F$.

\begin{lemma}\label{colagem de MGYM}
Let $D$ be an open subset of $\O$ with Lipschitz boundary, and let
$\mu$, $\nu \in L^\infty_w(\O \times Q;\M(\Rb^{d \times N}))$ be
such that $\{\mu_{(x,y)}\}_{(x,y) \in \O \times Q}$ and
$\{\nu_{(x,y)}\}_{(x,y) \in \O \times Q}$ are  two-scale gradient
Young measures with same underlying deformation $u \in
W^{1,p}(\O;\Rb^d)$. Let
$$\sigma_{(x,y)}:=\left\{\begin{array}{ll}
\mu_{(x,y)} & \text{if } ~ (x,y) \in D \times Q\\
\nu_{(x,y)} & \text{if } ~ (x,y) \in (\O \setminus D) \times Q.
\end{array}\right.$$
\noindent Then $\sigma \in L^\infty_w(\O \times Q;\M(\Rb^{d \times
N}))$ and $\{\sigma_{(x,y)}\}_{(x,y) \in \O \times Q}$ is  a
two-scale gradient Young measure with underlying deformation $u
\in W^{1,p}(\O;\Rb^d)$.
\end{lemma}
\begin{proof}
We have to show that for every sequence $\{\e_n\} \to 0$ there
exists $\{w_n\} \subset W^{1,p}(\O;\Rb^d)$ such that $w_n
\rightharpoonup u$ in $W^{1,p}(\O;\Rb^d)$ and   $\{(\la \cdot/\e_n
\ra, \nabla w_n)\}$ generates the Young measure $\{\sigma_{(x,y)}
\otimes dy\}_{x \in \O}$.

By Lemma \ref{boundary}, there exist sequences $\{u_n\} \subset
W^{1,p}(D;\Rb^d)$ and $\{v_n\} \subset W^{1,p}(\O \setminus
{\overline D};\Rb^d)$ such that $u_n \rightharpoonup u$ in
$W^{1,p}(D;\Rb^d)$, $v_n \rightharpoonup u$ in $W^{1,p}(\O \setminus
\overline D;\Rb^d)$, $u_n=v_n=u$ on $\partial D$ and such that
$\{(\la \cdot/\e_n \ra, \nabla u_n)\}$ and $\{(\la \cdot/\e_n \ra,
\nabla v_n)\}$ generate, respectively, the Young measures
$\{\mu_{(x,y)} \otimes dy\}_{x \in D}$ and $\{\nu_{(x,y)} \otimes
dy\}_{x \in \O\setminus \overline D}$.

Define
$$w_n:=\left\{\begin{array}{ll}
u_n & \text{if } x \in D,\\
v_n & \text{if } x \in \O \setminus \overline D
\end{array}\right.$$
Then $\{w_n\} \subset W^{1,p}(\O;\Rb^d)$, $w_n \rightharpoonup u$ in
$W^{1,p}(\O;\Rb^d)$ and given $z \in L^1(\O)$ and $\varphi \in
\C_0(\Rb^N \times \Rb^{d \times N})$ we have
\begin{eqnarray*}
&&\lim_{n \to +\infty} \int_{\O} z(x) \,
\varphi\left(\left\la\frac{x}{\e_n}\right\ra, \nabla w_n(x)\right)
\,dx\\
&&\hspace{2cm}= \lim_{n \to +\infty} \int_{D} z(x) \,
\varphi\left(\left\la\frac{x}{\e_n}\right\ra,
\nabla u_n(x)\right) \,dx \\
&&\hspace{2.5cm}+ \lim_{n \to +\infty} \int_{\O \setminus D} z(x)\,
\varphi\left(\left\la\frac{x}{\e_n}\right\ra,
\nabla v_n(x)\right) \,dx\\
&&\hspace{2cm}= \int_{\O} z(x) \int_Q \int_{\Rb^{d \times N}}
\varphi(x,\xi) \,d\sigma_{(x,y)}(\xi) \,dy\,dx,
\end{eqnarray*}
which concludes the proof.
\end{proof}

As a consequence of Remark \ref{indep} there is no loss of
generality to assume hereafter that $\O=Q$. We can now prove the
following result.

\begin{lemma}\label{closeconvex}
$M_{F}$ is a convex and weak*-closed subset of $(\E_p)'$.
\end{lemma}

\begin{proof}
We identify every element $\nu \in M_F$ with a homogeneous Young
measure $\nu_y \otimes dy$.

We start by showing that $M_{F}$ is a subset of $(\E_p)'$. For this
purpose let $\nu \in M_F$. Arguing exactly as
in the proof of Lemma \ref{nec} one can show that
$$K:=\int_Q \int_{\Rb^{d \times N}}|\xi|^p\,d\nu_y(\xi)\, dy <
+\infty.$$ Hence, using the fact that $\nu_y$ are probability
measures for a.e.\ $y \in Q$, for every $f \in \E_p$ we have that
\begin{eqnarray*}
\int_Q \int_{\Rb^{d \times N}}f(y,\xi)\, d\nu_y(\xi)\, dy & \leq &
\|f \|_{\E_p}\int_Q \int_{\Rb^{d \times
N}}(1+|\xi|^p)\,d\nu_y(\xi)\, dy\\
& = & (1+K)\|f\|_{\E_p}.
\end{eqnarray*}
As a consequence,  $M_F \subset (\E_p)'$.

Let us now prove that $M_F$ is closed for the weak*-topology of
$(\E_p)'$. Denoting by $\overline{M_F}$ the closure of $M_F$ for the
weak*-topology of $(\E_p)'$ it is enough  show that
$\overline{M_F}\subset  M_F$.  Since $\E_p$ is separable, the weak*-topology of $(\E_p)'$ is locally metrizable and
thus, if $\nu \in \overline{M_F}$, there exists a sequence
$\{\nu^k\} \subset M_F$ such that $\nu^k \xrightharpoonup{*}{} \nu$
in $(\E_p)'$. Hence, since the map $(y,\xi)\mapsto \xi_{ij}$ is in
 $\E_p$ (where $1\leq i \leq d$ and $1 \leq j \leq N$), we get,
from the definition of weak*-convergence in $(\E_p)'$, that
\begin{equation}\label{condF}
\int_Q \int_{\Rb^{d \times N}} \xi \,d\nu_y(\xi)\, dy=\lim_{k \to
+\infty}\int_Q \int_{\Rb^{d \times N}} \xi \,d\nu^k_y(\xi)\, dy=F.
\end{equation}
It remains to show that $\{{\nu_y}\}_{y\in Q}$ is a homogeneous
two-scale Young measure. By definition, given $\{\e_n\} \to 0$, for
each $k \in \Nb$ there exist sequences $\{u^k_n\}_{n \in \Nb}
\subset W^{1,p}(Q;\Rb^{d})$ such that $\{(\la\cdot/\e_n \ra, \nabla
u^k_n)\}_{n \in \Nb}$ generate the homogeneous Young measures
$\nu^k_y \otimes dy$. For every $(z,\varphi)$ in a countable dense
subset of $L^{1}(Q)\times \C_{0}(\Rb^N \times \Rb^{d\times N})$ we
have that

\begin{eqnarray*}&&\lim_{k\to +\infty}\lim_{n \to +\infty}
\int_{Q}z(x)\, \varphi\left(\left\la\frac{x}{\e_n}\right\ra, \nabla
u^k_n(x)\right)\, dx\\
&&\hspace{2cm}= \lim_{k \to +\infty}
\int_{Q}\int_{Q}\int_{\Rb^{d\times N}} z(x)\, \varphi(y,\xi)\,
d\nu^k_{y}(\xi) \,dy\, dx\\
&&\hspace{2cm}= \int_{Q}z(x)\, dx \int_{Q} \int_{\Rb^{d\times N}}
\varphi(y,\xi)\, d\nu_{y}(\xi) \,dy,
\end{eqnarray*}
where we have used the fact that $\C_0(\Rb^N \times \Rb^{d \times
N}) \subset \E_p$ in the second equality. By a diagonalization
argument we can find a sequence $\{k(n)\}\nearrow +\infty$ such
that, setting $v_n:=u^{k(n)}_n$, we have that
\begin{eqnarray*}
\lim_{n \to +\infty}\int_{Q}z(x)\,
\varphi\left(\left\la\frac{x}{\e_n}\right\ra, \nabla
v_n(x)\right)\, dx = \int_{\O}z(x)\, dx \int_{Q}\int_{\Rb^{d\times
N}} \varphi(y,\xi) \, d\nu_{y}(\xi) \,dy.
\end{eqnarray*}
Thus, $\{\nu_y\}_{y\in Q}$ is a homogeneous two-scale Young measure,
which together with (\ref{condF}) implies that $\nu\in M_F.$

Next we show that $M_{F}$ is convex. Given $\mu$, $\nu\in M_F$ and
$t \in (0,1)$ we have to show that  $t\mu+(1-t)\nu \in M_{F}$. Let
$D=(0,t)\times {(0,1)}^{N-1} \subset Q$ and define
$$\sigma_{(x,y)}:=\left\{\begin{array}{ll}
\mu_y & \text{if} ~ (x,y) \in D \times Q\\
\nu_y & \text{if} ~ (x,y) \in (Q \setminus D) \times Q.
\end{array}\right.$$
By Lemma \ref{colagem de MGYM} we have that
$\{\sigma_{(x,y)}\}_{(x,y) \in Q \times Q}$ is a two-scale
gradient Young measure and from Lemma \ref{average} its average
$\{\overline \sigma_y\}_{y \in Q}$ is a homogeneous two-scale
gradient Young measure. We claim that
$\overline{\sigma}=t\mu+(1-t)\nu$. Indeed, for every $\varphi \in
L^1(Q ; \C_0(\Rb^{d \times N}))$
\begin{eqnarray*}
\int_{Q}\int_{\Rb^{d \times N}}\varphi(y,\xi)\, d\overline
\sigma_y(\xi)\, dy &=& \int_{Q} \int_{Q}\int_{\Rb^{d \times
N}}\varphi(y,\xi)\,d\sigma_{(x,y)} (\xi)\, dy \,dx\\
& = & t\int_Q\int_{\Rb^{d \times N}}\varphi(y,\xi)\,
d\mu_y (\xi)\, dy\\
&& + (1-t) \int_Q\int_{\Rb^{d \times N}}\varphi(y,\xi)\, d\nu_y
(\xi)\, dy.
\end{eqnarray*}
In particular, $$\int_{Q}\int_{\Rb^{d \times N}}\xi\, d\overline
\sigma_y(\xi)\, dy =
t\int_Q\int_{\Rb^{d \times N}}\xi\, d\mu_y (\xi)\, dy\\
+ (1-t) \int_Q\int_{\Rb^{d \times N}}\xi\, d\nu_y (\xi)\, dy=F,$$
and thus $\overline{\sigma}=t\mu+(1-t)\nu \in M_F$.
\end{proof}

We are now in position to show the sufficiency of conditions i)-iii)
in (\ref{i})-(\ref{iii}) in the homogeneous case.

\begin{proof}[Proof of Lemma \ref{homo}]
Let $F\in \Rb^{d\times N}$ and $\nu\in L_{w}^{\infty}(Q;
\M(\Rb^{d\times N}))$ be such that $\nu_y \in {\mathcal
P}(\Rb^{d\times N})$ for a.e. $y\in Q$, and
(\ref{n1})-(\ref{powerp}) hold. We will proceed by contradiction
using the Hahn-Banach Separation Theorem. Assume that
$\{\nu_y\}_{y\in Q}$ is not a homogeneous two-scale Young measure.

By Lemma \ref{closeconvex}, $M_F$ is a convex and weak* closed
subset of $(\E_p)'$. Moreover, by (\ref{powerp}) and the fact that
$\{\nu_y\}_{y\in Q}$ is a family of probability measures, we get
that $\nu \in (\E_p)'$ as well (see {\it e.g.}\ the first part of
the proof of Lemma \ref{closeconvex}). As $\nu \not\in M_F$,
according to the Hahn-Banach Separation Theorem, we can separate
$\nu$ from $M_F$ \ie there exist a linear weak* continuous map $L :
(\E_p)' \to \Rb$ and $\alpha \in \Rb$ such that $\langle
L,\nu\rangle_{(\E_p)',\E_p} < \alpha$ and $\langle
L,\mu\rangle_{(\E_p)',\E_p} \geq \alpha$ for all $\mu \in M_F$. Let
$f \in \E_p$ be such that
\begin{eqnarray}\label{1}
\alpha \leq \langle L , \mu \rangle_{(\E_p)'',(\E_p)'} & = & \langle
\mu , f \rangle_{(\E_p)',\E_p}\nonumber\\
& = & \int_Q \int_{\Rb^{d \times N}} f(y,\xi)\, d\mu_y(\xi)\, dy
\quad\text{for all}\,\,\,\mu \in M_F,
\end{eqnarray}
and
\begin{eqnarray}\label{2}
\alpha >\langle L , \nu \rangle_{(\E_p)'',(\E_p)'} & = & \langle \nu
, f \rangle_{(\E_p)',\E_p}\nonumber\\
& = & \int_Q \int_{\Rb^{d \times N}} f(y,\xi)\, d\nu_y(\xi)\, dy
\geq f_{\rm hom}(F).
\end{eqnarray}
Let
$$f_{\rm H}(F) := \inf_{\mu \in M_F} \int_Q \int_{\Rb^{d \times N}}
f(y,\xi)\, d\mu_y(\xi)\, dy, \quad  F \in \Rb^{d \times N}.$$ Then,
by (\ref{1}), we have that  $\a \leq f_{\rm H}(F)$. We are going to
show that
\begin{equation}\label{n3}f_{\rm H}(F) \leq f_{\rm
hom}(F),\end{equation} which is a contradiction with (\ref{2}) and
asserts the conclusion of this lemma.

To prove (\ref{n3}),  let $T \in \Nb$ and $\phi \in
W_0^{1,p}((0,T)^N;\Rb^{d})$.  Extend $\phi$ to $\Rb^N$ by
$(0,T)^N$-periodicity and consider the sequence
$$\phi_n(x)= Fx+\e_n \phi\left(\frac{x}{\e_n}\right),$$
where $\{\e_n\} \to 0$ is an arbitrary sequence. Let $\varphi \in
\C_0(\Rb^N \times \Rb^{d \times N})$ and $z \in L^1(Q)$. Then, since
$T \in \Nb$, the function $y \mapsto \varphi(\la y \ra,F+\nabla
\phi(y))$ is $(0,T)^N$-periodic and according to the
Riemann-Lebesgue Lemma, we get that
\begin{eqnarray}
&&\lim_{n \to +\infty}\int_Q z(x)\,
\varphi\left(\left\la\frac{x}{\e_n}\right\ra,\nabla \phi_n(x)
\right)dx\nonumber\\
&&\hspace{2cm}= \lim_{n \to +\infty}\int_Q z(x)\,
\varphi\left(\left\la\frac{x}{\e_n}\right\ra,F+\nabla
\phi\left(\frac{x}{\e_n}\right) \right)dx\nonumber\\
&&\hspace{2cm}=\int_Q z(x)\,dx \med_{(0,T)^N} \varphi(\la y
\ra,F+\nabla \phi(y))\,dy.\label{n4}
\end{eqnarray}
Observe that
\begin{eqnarray}
&&\med_{(0,T)^N} \varphi(\la y \ra,F+\nabla
\phi(y))\,dy\nonumber\\
&&\hspace{1cm}= \frac1 {T^N} \sum_{a_i \in {\mathbb Z}^N \cap
[0,T)^N} \int_{a_i+Q} \varphi(\la y \ra ,F+\nabla
\phi(y))\,dy \nonumber\\
&&\hspace{1cm}= \frac1 {T^N} \sum_{a_i \in {\mathbb Z}^N \cap
[0,T)^N} \int_{Q} \varphi(\la a_i+y \ra,F+\nabla
\phi(a_i+y))\,dy \nonumber\\
&&\hspace{1cm}= \frac1 {T^N} \sum_{a_i \in {\mathbb Z}^N \cap
[0,T)^N} \int_{Q} \varphi(y,F+\nabla \phi(a_i+y))\,dy.\label{n5}
\end{eqnarray}
Thus, from (\ref{n4}) and (\ref{n5}), the pair $\{(\la\cdot
/\e_n\ra,\nabla \phi_n )\}$ generates the homogeneous Young measure
$$\mu:=\sum_{a_i \in {\mathbb Z}^N \cap [0,T)^N} \frac1 {T^N} \delta_{F+\nabla \phi(a_i+y)} \otimes dy.$$
Then
$$\int_Q \int_{\Rb^{d \times N}}\xi\, d\nu_{y}(\xi)\, dy = F,$$
which implies that $\mu\in M_F.$ In addition
$$\med_{(0,T)^N} f(\la y \ra,F+\nabla
\phi(y))\,dy=\int_{Q}\int_{\Rb^{d\times
N}}f(y,\xi)\,d\mu_{y}(\xi)\,dy,$$ and then
$$\med_{(0,T)^N}f(\la y\ra,F+\nabla
\phi(y))\,dy \geq f_{\rm H}(F).$$ As a consequence, taking the
infimum over all $\phi \in W_0^{1,p}((0,T)^N;\Rb^{d})$ and the limit
as $T \to +\infty$ we get that $f_{\rm hom}(F) \geq f_{\rm H}(F)$
which proves (\ref{n3}).
\end{proof}

Let us conclude this section by stating a localization result
which allows us to construct homogeneous two-scale gradient Young
measures starting from any kind of them.

\begin{proposition}\label{localization}
Let $\nu \in L^\infty_w(\O \times Q ; \M(\Rb^{d \times N}))$ be such
that $\{\nu_{(x,y)}\}_{(x,y) \in \O \times Q}$ is a two-scale
gradient Young measure. Then for a.e.\ $a \in \O$,
$\{\nu_{(a,y)}\}_{y \in Q}$ is a homogeneous two-scale gradient
Young measure.
\end{proposition}

\begin{proof}
Since $\{\nu_{(x,y)}\}_{(x,y) \in \O \times Q}$ is a two-scale
gradient Young measure, from Lemma \ref{nec} it satisfies
properties (\ref{i}), (\ref{ii}) and (\ref{iii}) above. Since
$u_1(x,\cdot)$ is $Q$-periodic for a.e.\ $x \in \O$, integrating
(\ref{i}) with respect to $y \in Q$, it follows that
\begin{equation}\label{ibis}
\int_Q \int_{\Rb^{d \times N}} \xi \, d\nu_{(x,y)}(\xi)\, dy= \nabla
u(x)\quad \text{ for a.e.\ }x \in \O.
\end{equation}
Furthermore, (\ref{iii}) implies that
\begin{equation}\label{iiibis}
\int_Q \int_{\Rb^{d \times N}} |\xi|^p\, d\nu_{(x,y)}(\xi)\, dy <
+\infty, \quad \text{ for a.e.\ }x \in \O.
\end{equation}
Let $E \subset \O$ be a set of Lebesgue measure zero  such that
(\ref{ibis}), (\ref{ii}) and (\ref{iiibis}) do not hold. Then for
every $a \in \O \setminus E$
$$\int_Q \int_{\Rb^{d \times N}} \xi \, d\nu_{(a,y)}(\xi)\, dy= \nabla
u(a),$$
$$\int_Q \int_{\Rb^{d \times N}} |\xi|^p \, d\nu_{(a,y)}(\xi)\, dy
<+\infty,$$ and
$$\int_Q \int_{\Rb^{d \times N}} f(y,\xi)\, d\nu_{(a,y)}(\xi)\, dy \geq
f_{\rm hom}( \nabla u(a))$$ \noindent for every $f \in \E_p$.

As a consequence of Lemma \ref{homo}, for every $a \in \O\setminus
E$, the family $\{\nu_{(a,y)}\}_{y \in Q}$ is a homogeneous
two-scale gradient Young measure.
\end{proof}

\subsubsection{The nonhomogeneous case}\label{nonhomogeneous}

\noindent We treat now the general case   whose proof is based on
Proposition \ref{localization} and a suitable decomposition of the
domain $\O$. We use (a variant of) Vitali's covering Theorem and an
approximation of two-scale gradient Young measures by measures of
this class that are piecewise constant   with respect to $x$.

\begin{lemma}\label{nonhomo}
Let $\O$ be a bounded and open subset of $\Rb^N$ with Lipschitz
boundary. Let $\nu \in L_w^\infty(\O \times Q;\M( \Rb^{d \times
N}))$ be such that $\nu_{(x,y)} \in \PP(\Rb^{d \times N})$ for
a.e.\ $(x,y) \in \O \times Q$. Suppose that
\begin{itemize}
\item[(i)] there exist $u \in W^{1,p}(\O;\Rb^d)$ and $u_1 \in
L^p(\O;W^{1,p}_{\rm per}(Q;\Rb^d))$ satisfying
\begin{equation}\label{1715}
\int_{\Rb^{d \times N}} \xi \, d\nu_{(x,y)}(\xi)= \nabla u(x) +
\nabla_y u_1(x,y)\quad \text{ for a.e.\ }(x,y) \in \O \times Q;
\end{equation}
\item[(ii)] for every $f \in \E_p$,
\begin{equation}\label{1721}
f_{\rm hom}(\nabla u(x)) \leq \int_Q \int_{\Rb^{d \times N}}
f(y,\xi)\, d\nu_{(x,y)}(\xi)\, dy \quad \text{ for a.e.\ }x \in \O;
\end{equation}
\item[(iii)] $\ds (x,y) \mapsto \int_{\Rb^{d \times N}} |\xi|^p\,
d\nu_{(x,y)}(\xi) \in L^1(\O \times Q).$
\end{itemize}
Then $\{\nu_{(x,y)}\}_{(x,y) \in \O \times Q}$ is a two-scale
gradient Young measure with underlying deformation $u$.
\end{lemma}

\begin{proof}
In a first step, we address the case where the underlying
deformation is zero,  while the general case is treated
afterwards.

{\it Step 1.} Assume $u=0$ and let $(\varphi, z)$ be in  a countable
dense subset of $\C_0(\Rb^N \times \Rb^{d \times N})\times L^1(\O)$.
Set
$$\overline \varphi(x):= \int_Q \int_{\Rb^{d \times N}}
\varphi(y,\xi)\, d\nu_{(x,y)}(\xi)\, dy.$$ Let $k \in \Nb$ and let
$E \subset \O$ be the set of Lebesgue measure zero given by
Proposition \ref{localization}. According to Lemma 7.9 in
Pedregal \cite{Pbook}, there exist points $a_i^k \in \O \setminus E$
and positive numbers $\rho_i^k \leq 1/k$ such that $\{a_i^k +
\rho_i^k \overline \O\}$ are pairwise disjoint for each $k$,
$$\overline \O=\bigcup_{i\geq 1}(a_i^k + \rho_i^k \overline \O) \cup
E_k,\qquad \LL^N(E_k)=0$$ and
\begin{equation}\label{Vitali}
\int_\O z(x) \, \overline \varphi(x)\, dx = \lim_{k \to +\infty}
\sum_{i\geq 1} \overline \varphi(a_i^k) \int_{a_i^k + \rho_i^k \O}
z(x)\, dx. \end{equation} For each $k \in \Nb$, let $m_k \in \Nb$
large enough so that
\begin{equation}\label{Mk}
\left|\sum_{i=1}^{m_k} \overline \varphi(a_i^k) \int_{a_i^k +
\rho_i^k \O} z(x)\, dx - \sum_{i \geq 1} \overline \varphi(a_i^k)
\int_{a_i^k + \rho_i^k \O} z(x)\, dx \right| < \frac{1}{k}.
\end{equation}
For fixed $i$ and $k$,  by the choice of  $a_i^k$ and  Proposition
\ref{localization} the family $\{\nu_{(a_i^k,y)}\}_{y \in Q}$ is a
homogeneous two-scale gradient Young measure. Hence by Remark
\ref{indep} and Lemma \ref{boundary}, for every sequence $\{\e_n\}
\to 0$, there exist sequences $\{u_n^{i,k}\}_{n \in \Nb} \subset
W_0^{1,p}(a_i^k + \rho_i^k\O;\Rb^d)$ such that
$$\lim_{n \to +\infty}\int_{a_i^k + \rho_i^k\O} z(x) \,
\varphi \left(\left\la \frac{x}{\e_n}\right\ra,\nabla u_n^{i,k}(x)
\right)\, dx = \overline \varphi (a_i^k) \int_{a_i^k + \rho_i^k\O}
z(x)\, dx .$$ Summing up
\begin{eqnarray}\label{homoaik}
&&\lim_{n \to +\infty}\sum_{i=1}^{m_k}\int_{a_i^k + \rho_i^k\O}
z(x)\, \varphi \left(\left\la \frac{x}{\e_n}\right\ra,\nabla
u_n^{i,k}(x) \right)\, dx\nonumber\\
&&\hspace{2cm} = \sum_{i=1}^{m_k}\overline \varphi (a_i^k)
\int_{a_i^k + \rho_i^k\O} z(x)\, dx.
\end{eqnarray}
Let us define
$$u_n^k(x):=\left\{
\begin{array}{ll}
u_n^{i,k}(x) & \text{ if } x \in a_i^k + \rho_i^k \O,\\
0 & \text{ otherwise}
\end{array}\right.$$
and remark that $u_n^k \in W^{1,p}_0(\O;\Rb^d)$. Since the sets
$a_i^k + \rho_i^k\O$ are pairwise disjoint for each $k$ we have
that
\begin{eqnarray}\label{compute} &&\int_\O z(x) \, \varphi
\left(\left\la \frac{x}{\e_n}\right\ra,\nabla u_n^{k}(x) \right)\,
dx\nonumber\\
&&\hspace{1cm}= \sum_{i \geq 1} \int_{a_i^k + \rho_i^k \O} z(x) \,
\varphi \left(\left\la \frac{x}{\e_n}\right\ra,\nabla u_n^{i,k}(x)
\right)\, dx\nonumber\\
&&\hspace{1cm}= \sum_{i = 1}^{m_k} \int_{a_i^k + \rho_i^k \O} z(x)\,
\varphi \left(\left\la \frac{x}{\e_n}\right\ra,\nabla
u_n^{i,k}(x) \right)\, dx\nonumber\\
&&\hspace{2cm}+ \int_{\O \cap \, \bigcup_{i > m_k} (a_i^k + \rho_i^k
\O)} z(x) \, \varphi \left(\left\la \frac{x}{\e_n}\right\ra,\nabla
u_n^{k}(x) \right)\, dx.
\end{eqnarray}
But as $z \in L^1(\O)$ and $\LL^N\left( \O \cap \,\bigcup_{i
> m_k}  (a_i^k + \rho_i^k \O) \right) \to 0$, as $k \to +\infty$, it
follows that
\begin{eqnarray}\label{meas}
\lim_{k \to +\infty}\lim_{n \to +\infty} \left| \int_{\O\cap\,
\bigcup_{i
> m_k} (a_i^k + \rho_i^k \O)} z(x)\, \varphi \left(\left\la
\frac{x}{\e_n}\right\ra,\nabla u_n^{k}(x) \right)\, dx \right|=0.
\end{eqnarray}
Then, gathering (\ref{Vitali})-(\ref{meas}) we obtain that
\begin{eqnarray*}
&&\lim_{k \to +\infty}\lim_{n \to +\infty} \int_\O z(x) \, \varphi
\left(\left\la \frac{x}{\e_n}\right\ra,\nabla u_n^{k}(x)
\right)\, dx\\
&&\hspace{1cm}  = \lim_{k \to +\infty}\lim_{n \to +\infty} \sum_{i =
1}^{m_k} \int_{a_i^k + \rho_i^k \O} z(x) \, \varphi \left(\left\la
\frac{x}{\e_n}\right\ra,\nabla u_n^{i,k}(x)
\right)\, dx\\
&&\hspace{1cm}  = \lim_{k \to +\infty} \sum_{i=1}^{m_k}\overline
\varphi (a_i^k) \int_{a_i^k + \rho_i^k\O} z(x)\, dx\\
&&\hspace{1cm}  = \lim_{k \to +\infty} \sum_{i\geq 1}\overline
\varphi (a_i^k) \int_{a_i^k + \rho_i^k\O} z(x)\, dx\\
&&\hspace{1cm}  =\int_\O z(x) \, \overline \varphi(x)\, dx.
\end{eqnarray*}
A diagonalization argument implies  the existence of a sequence
$\{k(n)\} \nearrow +\infty$, as $n \to +\infty$, such that upon
setting $u_n:=u_n^{k(n)}$, then
$$\lim_{n \to +\infty} \int_\O z(x) \, \varphi
\left(\left\la \frac{x}{\e_n}\right\ra,\nabla u_n(x) \right)\,
dx=\int_\O z(x)\, \overline \varphi(x)\, dx$$ and
$u_n\rightharpoonup 0$ in $W^{1,p}(\O;\Rb^d)$, which
completes the proof whenever $u=0$.\\

{\it Step 2.} Consider now a general $u \in W^{1,p}(\O;\Rb^d)$ and
$\nu$ satisfying properties (i)-(iii).   We define $\tilde \nu \in
L^\infty_w(\O \times Q; \M(\Rb^{d \times N}))$ by
\begin{equation}\label{nutilde}
\la \tilde \nu , \varphi \ra := \int_\O \int_Q \int_{\Rb^{d \times
N}} \varphi(x,y,\xi - \nabla u(x))\, d\nu_{(x,y)}(\xi)\, dy \, dx,
\end{equation} for every $\varphi \in L^1(\O \times Q;\C_0(\Rb^{d \times N}))$.
We can easily check that $\tilde \nu$ satisfies the analogue of
properties (i)-(iii) with $\tilde{u}=0$. Hence, applying Step 1,
for every sequence $\{\e_n\} \to 0$ we may find a sequence
$\{\tilde u_n\} \subset W^{1,p}(\O;\Rb^d)$ such that $\{(\la
\cdot/\e_n \ra,\nabla \tilde u_n)\}$ generates the Young measure
$\{\tilde\nu_{(x,y)} \otimes dy\}_{x \in \O}$. Defining
$u_n:=\tilde u_n + u$, we claim that $\{(\la \cdot/\e_n \ra,\nabla
u_n)\}$ generates $\{\nu_{(x,y)}\otimes dy\}_{x \in \O}$. Indeed
let $\psi \in L^1(\O;\C_0(\Rb^N \times \Rb^{d \times N}))$ and
define the $\tilde \psi(x,y,\xi):=\psi(x,y,\xi+\nabla u(x))$ where
$\tilde \psi \in L^1(\O;\C_0(\Rb^N \times \Rb^{d \times N}))$ as
well. Then by (\ref{nutilde}),
\begin{eqnarray*}
\lim_{n \to +\infty} \int_\O \psi\left(x,\left\la \frac{x}{\e_n}
\right\ra,\nabla u_n(x) \right)\, dx & = & \lim_{n \to +\infty}
\int_\O \tilde \psi\left(x,\left\la \frac{x}{\e_n}
\right\ra,\nabla \tilde u_n(x) \right)\, dx\\
& = & \int_\O \int_Q \int_{\Rb^{d \times N}}\tilde
\psi(x,y,\xi)\, d\tilde\nu_{(x,y)}(\xi)\, dy\, dx\\
& =& \int_\O \int_Q \int_{\Rb^{d \times N}} \psi(x,y,\xi)\,
d\nu_{(x,y)}(\xi)\, dy\, dx
\end{eqnarray*}
which completes the proof.
\end{proof}

The next corollary asserts the independence of the sequence in
Definition \ref{MGYM}.

\begin{coro}\label{scaleindep}
Let $\{u_n\}$ be a bounded sequence in $W^{1,p}(\O;\Rb^d)$. Assume
that there exists a sequence $\{\e_n\} \to 0$ such that the pair
$\{(\la \cdot/\e_n\ra,\nabla u_n)\}$ generates a Young measure
$\{\nu_{(x,y)} \otimes dy\}_{x \in \O}$. Then the family
$\{\nu_{(x,y)}\}_{(x,y) \in \O \times Q}$ is a two-scale gradient
Young measure.
\end{coro}

\section{Proof of Theorem \ref{exprGlim}}

\noindent Before proving Theorem \ref{exprGlim} we start by
recalling Valadier's notion of {\it admissible integrand}  (see
\cite{V}).

\begin{defi}\label{admint}
A function $f:\O\times Q \times \Rb^{d\times N}\rightarrow
[0,+\infty)$ is said to be an {\rm admissible integrand} if for any
$\eta>0$, there exist compact sets $K_{\eta}\subset \O$ and
$Y_{\eta}\subset Q$, with $\LL^{N}(\O\setminus K_{\eta})<\eta$ and
$\LL^{N}(Q\setminus Y_{\eta})<\eta$, and such that
$f|_{K_{\eta}\times Y_{\eta}\times \Rb^{d\times N}}$ is continuous.
\end{defi}

We observe that from Lemma 4.11 in Barchiesi \cite{Bar}, if $f$ is
an admissible integrand then, for fixed $\e>0$,   the function
$(x,\xi) \mapsto f(x,\la x/\e\ra ,\xi)$ is $\LL(\O) \otimes
\mathcal B(\Rb^{d \times N})$-measurable, where $\LL(\O)$ and
$\mathcal B(\Rb^{d \times N})$ denote, respectively, the
$\sigma$-algebra of Lebesgue measurable subsets of $\O$ and Borel
subsets of $\Rb^{d \times N}$. In particular, the functional
(\ref{main-funct}) is well defined in $W^{1,p}(\O;\Rb^{d})$.

\begin{proof}[Proof of Theorem \ref{exprGlim}]
Let $u\in W^{1,p}(\O;\Rb^{d})$ and let $\{\e_n\} \to 0$. We start
by showing that
\begin{equation}\label{glimsup}
\G\text{-}\limsup_{n \to +\infty}\F_{\e_n}(u)  \leq
\inf_{\nu \in \M_u} \int_\O \int_Q \int_{\Rb^{d \times N}}
f(x,y,\xi)\, d\nu_{(x,y)}(\xi)\,dy\, dx.
\end{equation}
where $\M_u$ is the set defined in (\ref{mu}). Let $\nu \in \M_u$,
by Remark \ref{underlying} there exists a sequence $\{u_n\}
\subset W^{1,p}(\O;\Rb^d)$ such that $\{(\la \cdot / \e_n \ra ,
\nabla u_n )\}$ generates the Young measure $\{\nu_{(x,y)} \otimes
dy\}_{x \in \O}$ and $u_n \rightharpoonup u$ in
$W^{1,p}(\O;\Rb^d)$. Extract a subsequence $\{\e_{n_k}\} \subset
\{\e_n\}$ such that
\begin{equation*}
\limsup_{n \to +\infty} \F_{\e_n}(u_n) = \lim_{k \to +\infty}
\F_{\e_{n_k}}(u_{n_k})
\end{equation*}
and that $\{|\nabla u_{n_k}|^p\}$ is equi-integrable, which is
always possible by the Decomposition Lemma (see Lemma 1.2 in
Fonseca, M\"uller \& Pedregal \cite{FMP}). In particular, due to
the $p$-growth condition (\ref{pgrowth}), the sequence
$\{f(\cdot,\la\cdot/\e_{n_k}\ra,\nabla u_{n_k})\}$ is
equi-integrable as well and applying Theorem 2.8 (ii) in Barchiesi
\cite{Bar2} we get that
\begin{eqnarray}\label{n7}
\G\text{-}\limsup_{n \to +\infty}\F_{\e_n}(u) & \leq & \lim_{k \to
+\infty} \int_\O
f\left(x,\left\la\frac{x}{\e_{n_k}}\right\ra,\nabla
u_{n_k}(x)\right)\, dx\\
& = &\int_\O \int_Q \int_{\Rb^{d \times N}} f(x,y,\xi)\,
d\nu_{(x,y)}(\xi)\, dy\, dx.
\end{eqnarray}
Taking the infimum over all $\nu \in \M_u$ in the right hand side
of the (\ref{n7}) yields to (\ref{glimsup}).

Let us prove now that
\begin{equation}\label{gliminf}
\G\text{-}\liminf_{n \to +\infty}\F_{\e_n}(u) \geq \inf_{\nu
\in \M_u} \int_\O \int_Q \int_{\Rb^{d \times N}} f(x,y,\xi)\,
d\nu_{(x,y)}(\xi)\,dy\, dx.
\end{equation}
 Let $\eta>0$ and $\{u_n\}\subset W^{1,p}(\O;\Rb^d)$ such that
$u_n \rightharpoonup u$ in $W^{1,p}(\O;\Rb^d)$ and
\begin{equation}\label{gliminfeta}
\liminf_{n \to +\infty} \F_{\e_n}(u_n) \leq \G\text{-}\liminf_{n
\to +\infty} \F_{\e_n}(u) + \eta.
\end{equation}
For a subsequence $\{n_k\}$, we can assume that there exists
$\nu \in L^\infty_w(\O \times Q;\M(\Rb^{d \times N}))$ such that
$\{(\la \cdot / \e_{n_k} \ra,\nabla u_{n_k})\}$ generates a
Young measure $\{\nu_{(x,y)} \otimes dy\}_{x \in \O}$ and
\begin{equation}\label{lim=liminf}
\lim_{k \to +\infty}\F_{\e_{n_k}}(u_{n_k}) = \liminf_{n \to
+\infty} \F_{\e_n}(u_n).
\end{equation}
We remark that $\{\nabla
u_{n_k}\}$ is equi-integrable since it is bounded in $L^p(\O;\Rb^{d
\times N})$ and $p>1$. Thus, by Theorem \ref{young} (v) we get that
for every $A \in \A(\O)$,
$$\int_A\nabla u(x)\, dx = \lim_{k \to +\infty} \int_A \nabla u_{n_k} (x)\,
dx= \int_A \int_Q \int_{\Rb^{d \times N}} \xi\, d\nu_{(x,y)}(\xi)\,
dy\, dx.$$ By the arbitrariness of the set $A$, it follows that
\begin{equation}\label{nablau(x)}
\nabla u(x) =\int_Q \int_{\Rb^{d \times N}} \xi\,
d\nu_{(x,y)}(\xi)\, dy\quad \text{ a.e. in }\O. \end{equation} As a
consequence of Corollary \ref{scaleindep} $\{\nu_{(x,y)}\}_{(x,y)
\in \O \times Q}$ is a two-scale gradient Young measure and, by
(\ref{nablau(x)}), we also have that $\nu \in \M_u$. Applying now
Theorem 2.8 (i) in Barchiesi \cite{Bar2} we get that
\begin{eqnarray*}
&&\lim_{k \to +\infty} \int_\O
f\left(x,\left\la\frac{x}{\e_{n_k}}\right\ra,\nabla u_{n_k}(x)
\right)\, dx\\
&&\hspace{2cm} \geq \int_\O \int_Q \int_{\Rb^{d \times N}} f(x,y,\xi)\, d\nu_{(x,y)}(\xi)\, dy\, dx\\
&&\hspace{2cm} \geq \inf_{\nu \in \M_u} \int_\O \int_Q \int_{\Rb^{d
\times N}} f(x,y,\xi)\, d\nu_{(x,y)}(\xi)\,dy\, dx.
\end{eqnarray*}
Hence by (\ref{gliminfeta}), (\ref{lim=liminf}) and the
arbitrariness of $\eta$ we get the desired result. Gathering
(\ref{glimsup}) and (\ref{gliminf}), we obtain that
$$\G\text{-}\lim_{n \to +\infty}\F_{\e_n}(u)=\inf_{\nu \in
\M_u} \int_\O \int_Q \int_{\Rb^{d \times N}} f(x,y,\xi)\,
d\nu_{(x,y)}(\xi)\,dy\, dx.$$

It remains to prove that the minimum is attained. To this aim,
consider a recovering sequence $\{\bar u_n\} \subset
W^{1,p}(\O;\Rb^d)$. Arguing exactly as before we can assume that (a
subsequence of) $\{\nabla \bar u_n\}$ generates a two-scale gradient
Young measure $\{\nu_{(x,y)}\}_{(x,y) \in \O \times Q}$, that $\nu
\in \M_u$ and $\{f(\cdot,\la\cdot/\e_n \ra,\nabla \bar u_n)\}$ is
equi-integrable. According to Theorem 2.8 (ii) in Barchiesi
\cite{Bar2} and using the fact that $\{\bar u_n\}$ is a recovering
sequence,
\begin{eqnarray*}
\G\text{-}\lim_{n \to +\infty}\F_{\e_n}(u) & = & \lim_{n \to
+\infty} \int_\O f\left(x,\left\la\frac{x}{\e_n}\right\ra,\nabla \bar u_n(x) \right)\, dx\\
& = & \int_\O \int_Q \int_{\Rb^{d \times N}} f(x,y,\xi)\,
d\nu_{(x,y)}(\xi)\,dy\, dx
\end{eqnarray*}
which completes the proof.
\end{proof}

Let us conclude by stating a Corollary which provides an alternative
formula to derive the homogenized energy density $f_{\rm hom}$ in
(\ref{fhom1}).

\begin{coro}\label{glim}
If $f:Q \times \Rb^{d \times N} \to [0,+\infty)$ is a Carath\'eodory
integrand (independent of $x$) and satisfying (\ref{pgrowth}), then
for every $u \in W^{1,p}(\O;\Rb^d)$,
$$\F_{\rm hom}(u)=\int_\O f_{\rm hom}(\nabla u(x))\, dx,$$
where for every $F \in \Rb^{d \times N}$,
$$f_{\rm hom}(F) = \min_{\nu \in M_F} \int_Q \int_{\Rb^{d \times N}} f(y,\xi)\, d\nu_y(\xi)\, dy$$
 and $M_F$ is defined in (\ref{MF}).
\end{coro}

\begin{proof}
It known from {\it e.g.}\ Theorem 14.5 in Braides \& Defranceschi
\cite{BD} that
$$\F_{\rm hom}(u)=\int_\O f_{\rm hom}(\nabla u(x))\, dx$$
where $f_{\rm hom}$ is defined in (\ref{fhom1}). By Theorem
\ref{exprGlim} with $\O=Q$ and $u(x)=Fx$, we get that
$$f_{\rm hom}(F)=\min_{\nu \in \M_u} \int_Q \int_Q \int_{\Rb^{d
\times N}} f(x,y,\xi)\, d\nu_{(x,y)}(\xi)\,dy\, dx.$$ The thesis
follows from Lemma \ref{average}.
\end{proof}

\vspace{1cm}

\noindent {\it Acknowledgments.} The authors wish to thank Irene
Fonseca for suggesting the problem. They also gratefully acknowledge
Gianni Dal Maso and Marco Barchiesi for fruitful comments and
stimulating discussions.

The research of J.-F. Babadjian has been supported by the Marie
Curie Research Training Network MRTN-CT-2004-505226 `Multi-scale
modelling and characterisation for phase transformations in advanced
materials' (MULTIMAT). He also acknowledge CEMAT, Mathematical Department of the Instituto Superior Técnico, Lisbon, for its hospitality and
support.

The research of M. Ba\'{\i}a was supported by Funda\c{c}\~{a}o
para a Ci\^{e}ncia e a Tecnologia under Grant PRAXIS XXI SFRH
$\hspace{-0.1cm}\backslash\hspace{-0.1cm}$ BPD
$\hspace{-0.1cm}\backslash\hspace{-0.1cm}$ 22775
$\hspace{-0.1cm}\backslash\hspace{-0.1cm}$ 2005 and by   Fundo
Social Europeu.


\begin{thebibliography}{99}





\bibitem{AM} {\sc G. Alberti \& S. M\"{u}ller}: A new approach to
variational problems with multiple scales, {\it Comm. Pure Appl.
Math.} {\bf 54} (2001), 761--825.

\bibitem{A} {\sc G. Allaire}: Homogenization and two-scale convergence, {\it SIAM J. Math.
Anal.} {\bf 23} (1992), 1482--1518.

\bibitem{AF} {\sc L. Ambrosio \& H. Frid}: Multiscale Young measures in almost periodic homogenization and
applications, Preprint CVGMT (2005).

\bibitem{HMM} {\sc O. Anza Hafsa, J.-P. Mandallena
\& G. Michaille}: Homogenization of periodic nonconvex integral
functionals in terms of Young measures, {\it ESAIM Control Optim.
Calc. Var.} {\bf 12} ( 2006), 35--51.

\bibitem{BF} {\sc M. Ba\'{\i}a \& I. Fonseca}: $\Gamma$-convergence of
functionals with periodic integrands via 2-scale convergence,
Preprint CNA (2005).

\bibitem{BF1} {\sc M. Ba\'{\i}a \& I. Fonseca}: The limit behavior of a
family of variational multiscale problems, to appear in {\it
Indiana. Univ. Math. J.}.

\bibitem{ball} {\sc J.M. Ball}: A version of the fundamental theorem for Young
measures, {\it PDE's and continuum models for phase transitions},
lecture notes physics {\bf 334} (1989), 207--215.

\bibitem{BJ} {\sc J.M. Ball \& R. James}: Fine phase mixtures as
minimizers of energy, {\it Arch. Rational Mech. Anal.} {\bf 100}
(1987), 15--52.

\bibitem{Bar} {\sc M. Barchiesi}: Multiscale homogenization of convex functionals
with discontinuous integrand, {\it J. Convex Anal.} {\bf 14} (2007),
205--226.

\bibitem{Bar2} {\sc M. Barchiesi}: Loss of polyconvexity by homogenization:
a new example, to appear in {\it Calc. Var. Partial Diff. Eq.}.

\bibitem{BC} {\sc E. Bonnetier \& C. Conca}: Approximation of Young measures by
functions and application to a problem of optimal design for plates
with variable thickness, {\it Proc. Roy. Soc. Edinburgh Sect. A.}
{\bf 124} (1994), 399--422.

\bibitem{B1} {\sc A. Braides}: Homogenization of some almost periodic coercive
functional, {\it Rend. Accad. Naz. Sci. XL.} {\bf 103} (1985),
313--322.

\bibitem{BLNotes} {\sc A. Braides:}  {\it A short introduction to Young
Measures}, Lecture notes SISSA (1999).

\bibitem{BD} {\sc A. Braides \& A. Defranceschi}: {\it Homogenization
of multiple integrals}, Oxford Lecture Series in Mathematics and its
Applications {\bf 12}, Oxford University Press, New York (1998).

\bibitem{BFF} {\sc A. Braides, I. Fonseca \& G. A. Francfort}: 3D-2D Asymptotic analysis for inhomogeneous thin films,
{\it Indiana Univ. Math. J.} {\bf 49} (2000), 1367--1404.

\bibitem{CK} {\sc M. Chipot \& D. Kinderlehrer}: Equilibrium
configurations of crystals, {\it Arch. Rational Mech. Anal.} {\bf
103} (1988), 237--277.

\bibitem{DP} {\sc R. J. DiPerna:} Convergence of approximate
solutions to conservation laws, {\it Arch. Rational Mech. Anal.}
{\bf 82} (1983), 27--70.

\bibitem{W} {\sc W. E:} Homogenization of linear and nonlinear
transport equations, {\it Comm. Pure Appl. Math.} {\bf 45} (1992),
301--326.

\bibitem{F} {\sc I. Fonseca}: Variational methods for elastic
crystals, {\it Arch. Rational Mech. Anal.} {\bf 97} (1987),
189--220.

\bibitem{FMP} {\sc I. Fonseca, S. M\"uller \& P. Pedregal}: Analysis of
concentration and oscillation effects generated by gradients, {\it
SIAM J. Math. Anal.} {\bf 29} (1998), 736--756.

\bibitem{KP2} {\sc D. Kinderlehrer \& P. Pedregal}: Characterization of Young measures
generated by gradients, {\it Arch. Rational Mech. Anal.} {\bf 115}
(1991), 329--365.

\bibitem{KP} {\sc D. Kinderlehrer \& P. Pedregal}: Gradient Young
measures generated by sequences in Sobolev spaces, {\it J. Geom.
Anal.} {\bf 4} (1994), 59--90.

\bibitem{LNW} {\sc D. Lukkassen, G. Nguetseng \& P. Wall}: Two scale
convergence, {\it Int. J. Pure Appl. Math.} {\bf 2} (2002), 35--86.

\bibitem{MP}  {\sc F. Maestre \& P. Pedregal}: Quasiconvexification in 3-D
for a variational reformulation of an optimal design problem in
conductivity, {\it Nonlinear Anal.} {\bf 64} (2006), 1962--1976.

\bibitem{MNRR} {\sc J. M\'alek, J. Ne\v cas, M. Rokyta \& M. R\r u\v zi\v cka}: {\it Weak and
measure-valued solutions to evolutionary PDEs}, Applied Mathematics
and Mathematical Compution {\bf 13}, Chapman \& Hall (1996).

\bibitem{Ma} {\sc P. Marcellini}: Periodic solutions and homogenization
of nonlinear variational problems, {\it Annali Mat. Pura Appl.} {\bf
117} (1978), 139--152.

\bibitem{Mu} {\sc S. M\"{u}ller}: Homogenization of nonconvex integral
functionals and cellular elastic materials, {\it Arch. Rational
Mech. Anal.} {\bf 99} (1987), 189--212.

\bibitem{M} {\sc S. M\"uller}: {\it Variational models for
microstructure and phase transitions}, Lecture Notes in Math. {\bf
1713} (1999), 85--210, Springer, Berlin.

\bibitem{N1} {\sc G. Nguetseng}: A general convergence result for a functional related to
 the theory of homogenization, {\it SIAM J. Math. Anal.} {\bf 20} (1989), 608--623.

\bibitem{Pbook} {\sc P. Pedregal}: {\it Parametrized measures and variational principles},
Progress in Nonlinear Differential Equations and their Applications
{\bf 30}, Birkhäuser Verlag, Basel (1997).

\bibitem{P2} {\sc P. Pedregal}: From micro to macroenergy through Young measures,
{\it Meccanica} {\bf 40} (2005), 329--338.

\bibitem{P3} {\sc P. Pedregal}: Vector variational problems and applications
to optimal design, {\it ESAIM Control Optim. Calc. Var.} {\bf 11}
(2005), 357--381.

\bibitem{P} {\sc P. Pedregal}: Multiscale Young measures,
{\it Trans. Am.  Math. Soc.} {\bf 358} (2006), 591--602.

\bibitem{TR} {\sc T. Roub\'{\i}\v{c}ek}: {\it Relaxation in optimization theory and variational
calculus}, Berlin, De Gruyter series in nonlinear analysis and
applications (1997)

\bibitem{T1} {\sc L. Tartar}: Compensated compactness and
applications to partial differential equations, {\it Nonlinear
analysis and Mechanics: Heriot-Watt Symposium}, Vol. IV (R.J Knops,
ed.), Reseach Notes in Math., Pitman, London  {\bf 39} (1979),
136--212.

\bibitem{V2} {\sc M. Valadier}: D\'esint\'egration d'une mesure sur un produit,
{\it C. R. Acad Ac. Paris} {\bf 276} S\'{e}rie A (1973), 33--35.

\bibitem{V1} {\sc M. Valadier}: Young measures, {\it Methods of
Nonconvex Analysis}, Lecture Notes in Math, Springer-Verlag, (1990),
152--188.

\bibitem{V} {\sc M. Valadier}: Admissible functions in two-scale
convergence, {\it Portugaliae Mathematica} {\bf 54} (1997),
148--164.

\bibitem{LCY3} {\sc L. C. Young}: {\it Lectures on the calculus of variations and optimal control theory},
W. B. Saunders, Philadelphia (1969).

\end{thebibliography}
\end{document}